\documentclass[preprint]{elsarticle}

\usepackage{amssymb}
\usepackage{amsmath}
\usepackage{amsthm}
\usepackage{mathpartir}
\usepackage[pagewise]{lineno}
\usepackage{graphicx}
\usepackage{nameref}
\newtheorem{Theorem}{Theorem}[section]

\newtheorem{proposition}[Theorem]{Proposition}
\newtheorem{lemma}[Theorem]{Lemma}

\newtheorem{corollary}[Theorem]{Corollary}
\theoremstyle{definition}
\newtheorem{definition}[Theorem]{Definition}
\newtheorem{example}[Theorem]{Example}

\newcommand{\proofofref}{}
\newproof{zproofof}{Proof of \proofofref}
\newenvironment{proofof}[1]
 {\renewcommand{\proofofref}{#1}\zproofof}
 {\endzproofof}

\def\AAA{\mathcal{A}}

\def\EEE{\mathcal{E}}

\def\HHH{\mathcal{H}}
\def\III{\mathcal{I}}
\def\JJJ{\mathcal{J}}

\def\LLL{\mathcal{L}}
\def\MMM{\mathcal{M}}
\def\NNN{\mathcal{N}}

\def\PPP{\mathcal{P}}

\def\epsilon{\varepsilon}
\def\phi{\varphi}
\def\rho{\varrho}
\def\theta{\vartheta}

\def\nat{\mathbb{N}}
\def\real{\mathbb{R}}
\def\bool{\mathbb{B}}
\def\iffdef{:\!\iff}
\def\imp{\Rightarrow}

\def\pot#1{\mathcal{P}(#1)}

\def\pmap{\stackrel{p}{\mapsto}}
\def\Pmap{\stackrel{p}{\longmapsto}}

\def\Qmap{\stackrel{q}{\longmapsto}}
\def\comp{\asymp}

\def\per{\approx}
\def\perp{\approx^{_{+}}}
\def\ptp{Pt}
\def\emb#1#2{emb_{#1}^{#2}}
\def\proj#1#2{proj_{#2}^{#1}}

\def\elem{\EEE\LLL}
\def\APX{\AAA\PPP}
\def\fil{\mathfrak{D}}

\def\limf{\underline{\textrm{lim}}}
\def\y#1{\mathsf{#1}}

\def\embb{\hookrightarrow}
\newcommand*\quot[2]{{{\textstyle #1}\big/{\textstyle #2}}}
\def\up{\,\uparrow\!}

\def\restr{\!\upharpoonright}

\def\y#1{\mathsf{#1}}

\journal{}

\begin{document}

\begin{frontmatter}

\title{Continuity in Potential Infinite Models}
\author{Matthias Eberl}
\ead{matthias.eberl@mail.de}

\affiliation{organization={LMU Munich},
            addressline={Theresienstr. 39}, 
            city={Munich},
            postcode={80333},
            country={Germany}}

\begin{abstract}
We introduce a model of simple type theory with potential infinite carrier sets. The functions in this model are automatically continuous, as defined in this paper. This notion of continuity does not rely on topological concepts, including domain theoretic concepts, which essentially use actual infinite sets. The model is based on the concept of a factor system, which generalizes direct and inverse systems. A factor system is basically an extensible set indexed over a directed set of stages, together with a predecessor relation between object states at different stages.  

The function space, when considered as a factor system, expands simultaneously with its elements. On the one hand, the space is subdivided more and more, on the other hand, the elements increase and are defined more and more precisely. In addition, a factor system allows the construction of limits that are part of its expansion process and not outside of it. At these limits, elements are indefinitely large or small, which is a contextual notion and a substitute for elements that are infinitely large or small (points). This dynamic and contextual view is consistent with an understanding of infinity as a potential infinite. 
\end{abstract}


\end{frontmatter}

\section{Introduction}

This paper is part of a project that aims to explain modern mathematics using the concept of a potential infinite. The potential infinite is viewed as an indefinitely extensible finite, making it a form of finitism. So the basic ``equation'' is this:
\[
\text{Potential infinite} = \text{Indefinite extensibility} + \text{Finitistic view}.
\]

It is worth noting that the finitistic perspective is not essential, as the approach can be easily generalized to other situations in which finite is replaced by a notion of being definite.\footnote{The finitistic perspective dates back to \cite{Mycielski1986} and \cite{lavine2009understanding}. Mycielski's approach is finitistic and does not explicitly mention the potential infinite. However, it is evident that his locally finite theories allow for a dynamic reading, which is required for the potential infinite. Lavine provided a philosophical explanation by interpreting the bounds, used in Mycielski's finitistic translation, as ``indefinitely large finite'' stages. There are other approaches to potential infinity, for example the one described in \cite{mostowski2003representing} differs from Mycielski's approach in that it is less dynamic. The approach presented in \cite{linnebo2019actual} uses modal logic instead.} Thus, the essential notion for formalizing the model is indefinite extensibility.

The basic aspect of the potential infinite is its dynamic character, which sees infinite totalities as processes rather than as completed sets. An infinite totality is not a completed, actual entity, but rather remains in a state of becoming. To formalize this idea we use generalizations of direct and inverse systems, called factor systems \cite{eberl2023}, which are closed under the function space construction. The potential infinite has a constructive aspect in a general sense, similar to the views in \cite{fletcher2007infinity} or \cite{schechter1996handbook}. However, as our approach is purely semantic, a constructive reading of the quantifier and logical connectives is not required, in particular, we do not reject classical axioms such as the excluded middle, axiom of choice, or any impredicative reasoning. Due to our finitistic view we go further than predicative approaches \cite{Feferman2007Predicativity}, which accepts countable totalities as completed.

Another important feature of the extensible structure is that both the elements and the space in which they exist expand simultaneously. The space is not established at the moment that all elements have been introduced, as in set theory, nor is there a completed space in which we then find or construct its elements, as in domain theory. The model thus differs from models such as the Kleene-Kreisel continuous functionals. An equivalent, seemingly more finitistic formulation, such as the hyperfinite type structure \cite{normann1999hyperfinite}, uses Fr\'{e}chet products, which also require actual infinite totalities.

The space in the extensible model presented here is divided into increasingly smaller areas while the elements are identified with increasing precision. This concept is similar to Brouwer's idea of (lawlike and lawless) choice sequences within a spread (\cite{brouwer1918begrundung}, \cite{brouwer1919begrundung}, the standard text on the subject is \cite{troelstra1979}), but without the aspect of constructive reasoning. Thus, the stages that we introduce do not refer to states of knowledge, as in Kripke models \cite{kripke1965semantical}, but to ``ontological'' states.

Indefinitely large stages in a factor system can be obtained through a process called \emph{compactification}, which involves a limit construction or, more generally, the introduction of a target. This limit extends the system in a structure-preserving way, while at the same time these limits are temporary and \emph{inside} the system. In other words, the introduction of limits is part of the indefinitely extensible character of the potential infinite. From the perspective of the elements, the limit elements are \emph{indefinitely} precisely determined but not infinitely precise, and from a spatial perspective, the areas are \emph{indefinitely} small rather than points, so we may see them as infinitesimals.\footnote{It is an infinitesimal in a relative sense. An infinitesimal is a number $a$ such that $|a| < \frac{1}{n}$ for all natural numbers. Since the universal quantifier is interpreted with a reflection principle \cite{eberl2022RML}, this condition says that $|a| < \frac{1}{n}$ for all $n \in \nat_i$ for some indefinitely large index $i$ relative to a context $C$ (formally written as $i \gg C$). The context $C$ depends on the state of the investigation, which includes the syntax as well as the model. For such an indefinitely large index $i$, the value of $|a|$ will be indefinitely small.}

Our approach therefore favors a \emph{creative process} over an \emph{all-embracing completeness} \cite{hellman2002maximality}, leading to a non-punctual understanding of the continuum with smallest extended parts. Being smallest is a relative notion, different, for example, from the mereological approach in \cite{hellman2018varieties}. It is crucial to avoid a na\"ive interpretation of the universal quantifier and instead apply a reflection principle \cite{eberl2022RML}.

In this paper we investigate the structure of these limits. A limit is naturally endowed with a family of partial equivalence relations (PERs) $\per_i$, approximating the identity relation. The equivalence classes of $\per_i$ then correspond to the increasingly smaller regions of the space. The definition of PERs on the function space thereby requires a new notion of continuity.

\subsection{A Motivating Example}

Let $\nat := \{0,1,2 \dots\}$, $\nat^+ := \nat \setminus \{0\}$ and $\nat_i := \{0, \dots, i-1\}$. For the sake of motivation, we start with the notion of uniform continuity of functions on real numbers. However, to make this an appropriate example, we will use binary representations $\sum_{j = 1}^{\infty} 2^{-j}b_j$ (with $b_j \in \{0,1\}$) of real numbers between $0$ and $1$. By abuse of language, let $[0,1)$ denote the set of these binary representations, ignoring the fact that two representations may refer to the same real number. 

A function $f : [0,1) \to [0,1)$ is uniform continuous iff (i.e., if and only if) for each real number $\epsilon > 0$ there is a real number $\delta > 0$ such that for all $a, b \in [0,1)$ with $\lvert a-b \rvert < \delta$ we have $\lvert f(a)-f(b) \rvert < \epsilon$. Note first that it is not necessary to take all real numbers $\epsilon > 0$ and $\delta > 0$, but only $\epsilon$ and $\delta$ from the set $\{2^{-i} \in \real \mid i \in \nat\}$. With $a \per_i b$ being an abbreviation for $\lvert a - b\rvert < 2^{-i}$, an equivalent definition is the following: $f : [0,1) \to [0,1)$ is uniformly continuous if for every $j \in \nat$  there is a  $i \in \nat$ such that for all $a, b \in [0,1)$, $a \per_i b$ implies $f(a) \per_j f(b)$.

We do not need such a number $i \in \nat$ for every $j \in \nat$, but only for all $j \in \NNN$ in a cofinal subset $\NNN \subseteq \nat$. Cofinality of $\NNN$ means that $\forall j \in \nat \, \exists j' \geq j$ with $j' \in \NNN$. Moreover, it is obvious that if we have found such a number $i$ for some $j$, then every $i' \geq i$ also satisfies the condition. So the uniform continuity condition can be roughly described as follows: We must have \emph{sufficiently many} pairs $(i,j) \in \nat \times \nat$, such that for all $a, b \in \real$, $a \per_i b$ implies $f(a) \per_j f(b)$. In Section \ref{tarsec} we introduce such a notion of ``sufficiently many'' or ``indefinitely many'', which allows us to deal with function spaces in general.

In order to motivate the concepts of Section \ref{peruniv}, we will now switch to another, narrower definition of $\per_i$. Let $\MMM_i$ be the set of the $2^i$ intervals $[2^{-i}n,2^{-i}(n+1))$ for $n \in \nat_{2^i}$ and let $a \per_i b$ hold in the case that $a,b \in A_i$ for some interval $A_i \in \MMM_i$. For $i' \geq i$ let $\proj{i'}{i} : \MMM_{i'} \to \MMM_i$ map an interval $A_{i'} \in \MMM_{i'}$ to the unique interval $A_i \in \MMM_i$ in which it lies, i.e., which satisfies $A_{i'} \subseteq A_i$. $(\MMM_i)_{i \in \nat}$ with the projections $\proj{i'}{i'}$ defines an inverse system of nested intervals. The inverse limit consists of the binary representations of $a \in [0,1)$ together with the projections $Proj_i : [0,1) \to \MMM_i$ of the limit to the $i$-th approximation, defined by $Proj_i(a) = A_i \iff a \in A_i$.

In addition to the intervals, we can introduce ``grid points'' on $[0,1)$. These are the points $\PPP_i := \{2^{-i}n \in [0,1) \mid n \in \nat_{2^i}\}$ at stage $i$, which is the leftmost point of each interval. For all $a \in [0,1)$, there exists its approximation $pt_i(a) \in \PPP_i$ at stage $i$, which is the unique point $a \in \PPP_i$ with $pt_i(a) \per_i a$. If $a = \sum_{j = 1}^{\infty} 2^{-j}b_j$, then $pt_i(a)$ is the partial sum $\sum_{j = 1}^{i} 2^{-j}b_j$. We can therefore identify $[0,1)$ with the set of all of these convergent sequences $pt_i(a)_{i \in \nat}$ and equally take them instead of the inverse system of nested intervals. There are also embeddings $Emb_i : \MMM_i \to [0,1)$ defined by $A_i \mapsto 2^{-i}n$ for $A_i = [2^{-i}n,2^{-i}(n+1))$, which map to an interval its leftmost point. Both maps form an embedding-projection pair, satisfying $Proj_i \circ Emb_i = id_{\MMM_i}$ and $Emb_i \circ Proj_i = pt_i$.

\subsection{Some Observations}
\label{observesec}

In the example above, the space is described by the inverse system of nested intervals, while the elements in this space are the convergent sequence of partial sums. The dynamic aspect of the real line and its nested intervals is its increasing divisibility into smaller and smaller intervals. The dynamic aspect of the points on the number line is their determination with increasing precision. This is also an inverse system, isomorphic to that of the nested intervals.

Both processes, bisecting intervals and increasing the precision by one bit, are isomorphic as factor systems. Therefore, they have the same limit structure, which we identify. This identification corresponds to the common idea that an infinite process of nested intervals and the infinite process of real number representation lead to the same real number in the limit. In our approach, the limits of both processes are the same as well, but this limit consists of \emph{indefinitely} fine divided intervals from the perspective of space and \emph{indefinitely} exact determined values of real numbers from the perspective of their elements. If we regard the limit intervals as infinitesimals, then different to non-standard analysis, no new numbers for these infinitesimals are required.\footnote{For an overview of infinitesimals, see \cite{ApproachestoAnalysisInfinitesimals}.} 

Observe that it is not necessary to include all possible divisions of $[0,1)$ and all degrees of precision $2^{-i}$ to establish continuity. It is sufficient to meet the requirements of the definition of continuity, as determined by the indices $i,j \in \nat$ that specify the precision of the domain and codomain. Since we assume that the elements and the space, in which they exist, expand in a simultaneous process (as partial sums and as nested intervals), there is no need to have infinitely many divisions of the real number line either.

\subsection{Structure of the Paper}

Section \ref{factsec} briefly reviews the definition and main properties of a factor system. An important notion is that of a target and a limit of a factor system, introduced in Section \ref{tarsec}. A limit is a target satisfying a universal property. A target is defined with respect to a factor system and inherits the structure of the underlying factor system. Furthermore, the introduction of a target can be seen as a temporary extension of the factor system, not as a final state of completion. This change of perspective is one of the key ideas for the consistent realization of the potential infinite. 

There are some more general, afferent concepts, namely prefactor systems and prefactor systems with embeddings/projections, which are described in the beginning. These structures are however not closed under the function space construction.

Section \ref{peruniv} is the main part of this paper and deals with the structure of limits of factor systems independent of the underlying factor system. We gradually introduce more structural elements and properties on these limits. In Section \ref{prelimsec} we consider families of PERs (partial equivalence relations) that satisfy basic properties of an equality, called $\fil$-sets. In a second step we introduce points on these PERs (Section \ref{richuni}) and finally extensions of PERs to (total) equivalence relations in Section \ref{eqrelpersec}. 

\begin{figure}[ht]
\centering
\includegraphics[trim = 0mm 0mm 0mm 0mm, width=1\textwidth]{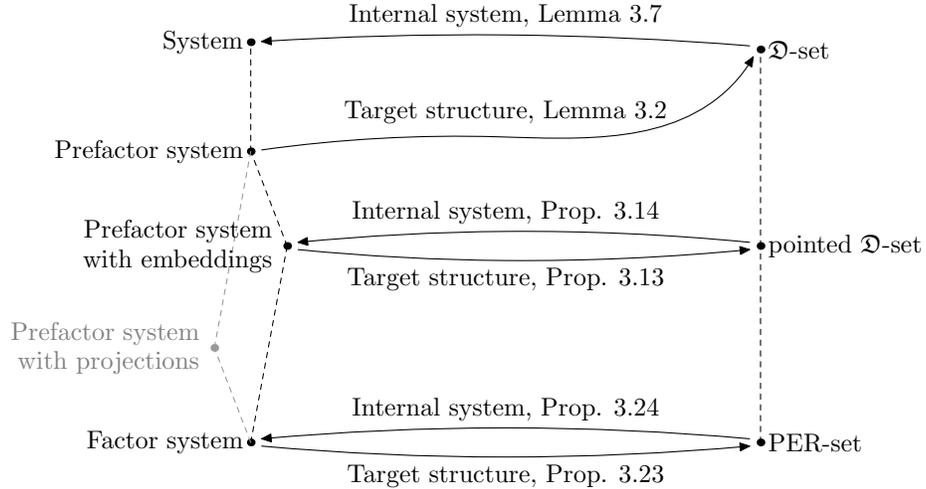}
\caption{Overview over the structures.}
\label{fig0}
\end{figure}

A \emph{PER-set} then has all of these additional structural elements, which are needed in order to carry out the function space construction in Section \ref{funPERset}. We also consider the reverse construction, inverse to the introduction of a target for a factor system, giving the \emph{internal system} of a PER-set, $\fil$-set or a pointed $\fil$-set respectively. Figure \ref{fig0} gives an overview of the connections between the systems on the one hand and the structure of their targets on the other hand.

\section{Factor Systems and their Limits}
\label{factlimsec}

Factor systems and their limits were introduced in \cite{eberl2023}. Here we recapitulate their definition and properties, but do not prove anything new.

\subsection{Factor Systems}
\label{factsec}

$\III$, $\JJJ$, $\HHH$ will always denote a non-empty directed index set with preorder $\leq$. By a \emph{system} we mean a pair $(\MMM_\III,\pmap)$ consisting of a family $\MMM_\III := (\MMM_i)_{i \in \III}$, finite sets $\MMM_i$ and reflexive (for $i = i'$) relations $\pmap$ on $\MMM_{i'} \times \MMM_i$ for $i' \geq i$. Two elements $a_i \in \MMM_{i}$ and $b_j \in \MMM_{j}$ are \emph{consistent}, written as $a_i \comp b_j$, iff there is an index $i' \geq i, j$ and some $a_{i'} \in \MMM_{i'}$ such that $a_{i'} \pmap a_i$ and $a_{i'} \pmap b_j$. As a convention, whenever we use a suffix $i \in \III$ for some element, this refers to the state, e.g.~$a_i \in \MMM_i$. An important special case is that the relations $\pmap$ are partial functions, which is equivalent to:
\begin{align}
\tag{Fun}
\label{idmapeq}
&a_i \comp b_i \iff a_i \pmap b_i \iff a_i = b_i
\end{align}
for all $a_i, b_i \in \MMM_i$ and all $i \in \III$. A system that satisfies (\ref{idmapeq}) is called \emph{standard}. $(\MMM_\III,\pmap)$ is a \emph{prefactor system} iff it is a system that satisfies
\begin{equation}
\label{pmapeq} 
\tag{Factor} 
a_{i'} \comp a_i \iff a_{i'} \pmap a_i
\end{equation}
for all $a_i \in \MMM_{i}$ and $a_{i'} \in \MMM_{i'}$ with $i \leq i'$. The relation $\comp$ is then an equivalence relation on a single set $\MMM_i$ with $a_i \comp b_i \iff a_i \pmap b_i \iff b_i \pmap a_i \ \text{ for } a_i, b_i \in \MMM_i$. In a prefactor system $b_{i'} \comp a_{i'} \pmap a_i$ implies $b_{i'} \pmap a_i$. A system $(\MMM_\III,\pmap)$ is called \emph{stable} iff for all $i' \geq i$, all $a_{i'} \in \MMM_{i'}$ and $a_i, b_i \in \MMM_i$
\begin{equation}
\label{stabeq} 
\tag{Stab} 
a_{i'} \pmap a_i \comp b_i \ \imp \ a_{i'} \pmap b_i.
\end{equation}

A family $emb = (\emb{i}{i'})_{i \leq i'}$ of $\comp$-embeddings consists of $\comp$-preserving maps $\emb{i}{i'} : \MMM_i \to \MMM_{i'}$ satisfying $\emb{i}{i}(a_i) \comp a_i$ and $\emb{i'}{i''}(\emb{i}{i'}(a_i)) \comp \emb{i}{i''}(a_i)$. The requirement $\comp$-preserving means that $\emb{i}{i'}(a_i) \comp \emb{i}{i'}(b_i)$ holds for $a_i \comp b_i$, for all $a_i, b_i \in \MMM_i$. A similar definition applies to $\comp$-projections $proj = (\proj{i'}{i})_{i \leq i'}$. Moreover, $\comp$-embeddings $emb$ together with $\comp$-projections $proj$ form an \emph{$\comp$-embedding-projection pair} iff $\proj{i'}{i}(\emb{i}{i'}(a_i)) \comp a_i$ holds for all $a_i \in \MMM_i$ and all $i \leq i'$. 

The $\comp$-embeddings $emb$ and $\comp$-projections $proj$ are \emph{coherent} if they satisfy for all $i \leq i' \leq i''$
\begin{align}
\label{embcond}
\tag{Emb}
a_{i'} \pmap a_i \ &\imp\ \emb{i'}{i''}(a_{i'}) \pmap a_i \text{ and} \\
\tag{Proj}
\label{projcond}
a_{i''} \pmap a_i \ &\imp\ \proj{i''}{i'}(a_{i''}) \pmap a_i \text{ resp.}
\end{align}

A $\comp$-embedding-projection pair $(emb,proj)$ is coherent if $emb$ and $proj$ are both coherent. Property (\ref{embcond}) implies that $\emb{i}{i'}(a_i) \pmap a_{i}$ holds for all $a_i \in \MMM_i$, and in case that Property (\ref{idmapeq}) holds, $\pmap$ is then a \emph{partial surjection}.

A \emph{prefactor system with embeddings} is a prefactor system $(\MMM_\III,\pmap)$ with coherent $\comp$-embeddings $emb$. A \emph{prefactor system with projections} is a prefactor system $(\MMM_\III,\pmap)$ with coherent $\comp$-projections $proj$.  A \emph{factor system} is a prefactor system $(\MMM_\III,\pmap)$ with a coherent $\comp$-embedding-projection pair $(emb,proj)$.

\begin{example}
\label{natex}
Consider $(\nat_i)_{i \in \nat^+}$ with the embedding-projection pair $\emb{i}{i'} : \nat_i \to \nat_{i'}$, $n \mapsto n$ and $\proj{i'}{i} : \nat_{i'} \to \nat_i$ with $n \mapsto \min(n,i-1)$. Then $((\nat_i)_{i \in \nat^+},\pmap,emb,proj)$, with $\nat_{i'} \ni n \pmap n \in \nat_i$ for all $n < i$, is a factor system.
\end{example}

A \emph{homomorphism} $\Phi = (\Phi_0, (\Phi^i)_{i \in \III})$ between two systems $(\MMM_\III, \pmap)$ and $(\NNN_\JJJ, \pmap)$ consists of maps $\Phi_0 : \III \to \JJJ$ and $\Phi^i : \MMM_i \to \NNN_{\Phi_0(i)}$ such that $\Phi_0$ is monotone and for all $i \leq i'$
\begin{equation}
\label{homeq}
a_{i'} \pmap a_i \ \imp \ \Phi^{i'}(a_{i'}) \pmap \Phi^i(a_i).
\end{equation}

If the equivalence holds in (\ref{homeq}), then $\Phi$ is said to be \emph{strong}. A homomorphism between two prefactor systems with embeddings $(\MMM_\III, \pmap, emb)$ additionally satisfies $\Phi^{i'}(\emb{i}{i'}(a_i)) = \emb{j}{j'}(\Phi^i(a_i))$, $j := \Phi_0(i)$, $j' := \Phi_0(i')$, a homomorphism between two prefactor systems with projections $(\MMM_\III, \pmap, proj)$ satisfies $\Phi^{i}(\proj{i'}{i}(a_{i'})) = \proj{j'}{j}(\Phi^{i'}(a_{i'}))$, and a homomorphism between two factor systems satisfies both equations. We call $\Phi$ injective (surjective) iff every map in $\Phi$ is injective (surjective resp.). An \emph{isomorphism} $\Phi$ between two systems (prefactor systems with embeddings/projections, factor systems) is a homomorphism with a further homomorphism $\Psi$ as its inverse, i.e., each part of $\Psi$ is inverse to that of $\Phi$. We write $\MMM_\III \simeq \NNN_\JJJ$ if an isomorphism between $\MMM_\III$ and $\NNN_\JJJ$ exists.

The function space of two factor systems $\MMM_\III$ and $\NNN_\JJJ$, denoted as $[\MMM_\III \to \NNN_\JJJ]$, is a family of (finite) sets $[\MMM_i \to \NNN_j]$ indexed by pairs $(i,j) \in \III \times \JJJ$ with product order, whereby we write $i \to j$ for such an index in $\III \times \JJJ$. The set $[\MMM_i \to \NNN_j]$ consists of all (total) functions $f : \MMM_i \to \NNN_j$ which preserve $\comp$, i.e., which satisfy $f(a_i) \comp f(b_i)$ for $a_i \comp b_i$. If the relations $\pmap$ are partial functions on $\MMM_\III$ and $\NNN_\JJJ$, then $\comp$ is the identity on $\MMM_i$ and $\NNN_j$ and $[\MMM_i \to \NNN_j]$ simply consists of all functions from $\MMM_i$ to $\NNN_j$.

Let $f \in [\MMM_{i} \to \NNN_{j}]$, $f' \in [\MMM_{i'} \to \NNN_{j'}]$ and $i \to j \leq i' \to j'$, i.e., $i \leq i'$ and $j \leq j'$. The basic relation $\pmap$ on the function space is a logical relation, c.f.~\cite{plotkin1973lambda}. It is thus defined as
\begin{equation*}
f' \pmap f \ \iffdef \ a_{i'} \pmap a_i \text{ implies } f'(a_{i'}) \pmap f(a_i)
\end{equation*}
for all $a_{i'} \in \MMM_{i'}$ and $a_i \in \MMM_{i}$. The embedding-projection pair for the function space is defined in the usual way:
\begin{align*}
\emb{i \to j}{i' \to j'} : [\MMM_i \to \NNN_j] \to [\MMM_{i'} \to \NNN_{j'}] &\quad f \mapsto \emb{j}{j'} \circ f \circ \proj{i'}{i},\\
\proj{i' \to j'}{i \to j} : [\MMM_{i'} \to \NNN_{j'}] \to [\MMM_i \to \NNN_j] &\quad f' \mapsto \proj{j'}{j} \circ f' \circ \emb{i}{i'}.
\end{align*}

If $\MMM_\III$ and $\NNN_\JJJ$ are both factor systems, so is their function space. Moreover, if $\pmap$ satisfies (\ref{idmapeq}) or (\ref{stabeq}) on $\NNN_\JJJ$, so does $\pmap$ on $[\MMM_\III \to \NNN_\JJJ]$.

\subsection{Targets and Limits}
\label{tarsec}

An interpretation of $\lambda$-terms in a potential infinite model is primarily an interpretation in a factor system. However, to prove the consistency of this interpretation --- that the value at two different stages $i$ and $j$ is consistent in the sense defined in Section \ref{factsec} --- a target of the factor system is required. The existence of a stage $i' \geq i, j$ is ensured by the directedness of the index set, and the concept of indefinitely many indices (used in the definition of a target) in addition guarantees the existence of a stage $i'$ where two considered elements are defined.

The subsequent concepts require a set $\fil(\III) \subseteq \pot{\III}$ --- $\pot{\III}$ denotes the power set of $\III$. $\HHH \in \fil(\III)$ states that there are \emph{indefinitely many}, or \emph{sufficiently many} indices in $\HHH$. The locution ``$\fil$-many indices'' refers to a set in $\fil(\III)$, and ``cofinal many indices'' refers to a cofinal index set. This notion applies to products $\III_0 \times \dots \times \III_{n-1}$ with the product order as well, giving the sets $\fil(\III_0 \times \dots \times \III_{n-1})$. We require that the notion of indefinitely many is in general stronger than cofinality, that is,
\begin{equation}
\label{cofineq}
\tag{Cofin}
\fil(\III) \subseteq \{\HHH \subseteq \III \mid \HHH \text{ is cofinal}\},
\end{equation}
including the situation that $\III$ is $\III_0 \times \dots \times \III_{n-1}$. A \emph{target} $(\MMM,\Pmap)$ for a system $\MMM_\III$ extends the system ``at the top'', i.e., the extension leads to a system $\MMM_{\bar\III}$, called \emph{compactification} of $\MMM_\III$, with index set $\bar\III := \III \cup \{top\}$, $top$ as greatest index, and $\MMM_{top} = \MMM$. Sometimes we will also refer to a target as a $\fil$-target in order to emphasize the dependency on $\fil(\III)$. 

Let $a \Pmap a_i$ denote $a \pmap a_i$ provided that $a \in \MMM$, and we also write $Emb_i$ for $\emb{i}{top}$ and $Proj_i$ for $\proj{top}{i}$. Relation $\pmap$ on a target $\MMM$, and consequently relation $\comp$ on $\MMM$ as well, is by definition the identity. The \emph{extension} of an element $a \in \MMM$ is $Ext(a) := \{a_i \in \bigcup_{i \in \III} \MMM_i \mid a \Pmap a_i\}$. A target $\MMM$ for a system $\MMM_\III$ has to satisfy 
\[
\III_a := \{i \in \III \mid \exists a_i \in \MMM_i \ a \Pmap a_i\} \in \fil(\III)
\]
for all objects $a \in \MMM$. Moreover, if the system is a prefactor system, a prefactor system with embeddings/projections or a factor system, then the compactification $\MMM_{\bar\III}$ must have this additional structure with its properties as well. Whereas the compactification of a system is automatically a system, the compactification of a \emph{prefactor system} requires for all $a \in \MMM$, $a_i \in \MMM_i$, $a_{i'} \in \MMM_{i'}$ and $i \leq i'$
\begin{equation}
\label{ppmapeq}
a \Pmap a_{i'}, \, a \Pmap a_i \ \imp \ a_{i'} \pmap a_i.
\end{equation}

If $\MMM$ is a target for a \emph{prefactor system with embeddings} $\MMM_{\III}$, then this implies the existence of $\comp$-embeddings $Emb_i : \MMM_i \to \MMM$, satisfying $Emb_{i'}(\emb{i}{i'}(a_i)) = Emb_{i}(a_i)$ and 
\begin{equation}
\label{embtareq}
a_{i'} \pmap a_i \ \imp\ Emb_{i'}(a_{i'}) \Pmap a_i \text{ for all } i \leq i',
\end{equation}
$a_i \in \MMM_i$, $a_{i'} \in \MMM_{i'}$ and $a \in \MMM$. If $\MMM$ is a target for a \emph{prefactor system with projections} $\MMM_{\III}$, then there are moreover $\comp$-projections $Proj_i : \MMM \to \MMM_i$ such that  $\proj{i'}{i}(Proj_{i'}(a)) \comp Proj_i(a)$ and $a \Pmap a_i \ \imp\ Proj_{i'}(a) \pmap a_i$. If $\MMM$ is a target for a \emph{factor system}, then $Emb$ and $Proj$ with these properties exist and both form a $\comp$-embedding-projection pair. We write $(\MMM,\Pmap)$, $(\MMM,\Pmap,Emb)$, $(\MMM,\Pmap,Proj)$ and $(\MMM,\Pmap,Emb,Proj)$, respectively, for these targets. 

Let $\MMM$ be a target for a prefactor system with projections $\MMM_\III$, then the projections $Proj_i$ can be defined by 
\begin{align}
\label{Projdef}
Proj_i(a) &:= \proj{i'}{i}(a_{i'}) \text{ for some } i' \geq i \text{ with } a \Pmap a_{i'}
\end{align}
for $a \in \MMM$. The properties of a prefactor system guarantee that $Proj_i(a)$ is uniquely defined modulo $\comp$. Since $Emb_{i'}(a_{i'}) \Pmap a_{i'}$ for $a_{i'} \in \MMM_{i'}$ we have for $i \leq i'$
\begin{equation}
\label{Projprop}
Proj_i(Emb_{i'}(a_{i'})) \comp \proj{i'}{i}(a_{i'}).
\end{equation}

A target $(\MMM,\Pmap)$ for a system $\MMM_\III$ is a \emph{limit} (or $\fil$-limit) of $\MMM_\III$ iff for every further target $(\NNN,\Qmap)$ for $\MMM_\III$ there is a unique map $\Phi : \NNN \to \MMM$ such that $a \Qmap a_i$ implies $\Phi(a) \Pmap a_i$. If the underlying system is a factor system or a prefactor system, then we call the limit \emph{factor limit} and \emph{prefactor limit}, resp. It turns out, however, that a factor limit is the same as the prefactor limit, and that this limit $\limf(\MMM_\III)$ is unique modulo isomorphism. Therefore we simply speak of a \emph{limit}, or a \emph{limit set}, if we want to distinguish it from a \emph{limit element} in this limit set.

It is possible to define concrete targets and limits in the form of sets of states: A set $\alpha \subseteq \bigcup_{i \in \III} \MMM_i$ in a system $\MMM_\III$ is called a \emph{consistent set} iff $a_{i'} \pmap a_i$ holds for all $a_{i'}, a_i \in \alpha$ with $i' \geq i$, and $\III_\alpha := \{i \in \III \mid \alpha \cap \MMM_i \not= \emptyset\} \in \fil(\III)$. If we already have a target $\MMM$ for a prefactor system, then the set $Ext(a)$ is such a consistent set for all $a \in \MMM$. The set of all consistent sets in a (pre)factor system $\MMM_\III$ is itself a target for $\MMM_\III$ with $\alpha \Pmap a_i \iffdef a_i \in \alpha$. In other words, if $\MMM$ denotes the set of all consistent sets, then the target is $(\MMM, \ni)$. 

A \emph{dynamic element} is a maximal (w.r.t.~subset inclusion) consistent set, and for each consistent set $\alpha$ in a prefactor system there is exactly one dynamic element $\alpha^m$ such that $\alpha \subseteq \alpha^m$. Whenever $\alpha^m = \beta^m$ holds for two consistent sets $\alpha$ and $\beta$, we write $\alpha \sim \beta$. Let $\elem(\MMM_\III)$ denote the set of all of these dynamic elements, then 
\begin{equation}
\label{elemtergeteq}
(\elem(\MMM_\III),\ni)
\end{equation}
is a (pre)factor limit of $\MMM_\III$.

Call an element $a \in \MMM$ of a target $\MMM$ for a prefactor system $\MMM_\III$ a \emph{limit (element)} of a consistent set $\alpha$ iff $\alpha \sim Ext(a)$. Then $a$ is a limit of its extension $Ext(a)$. If there is only one limit element, we denote it as $lim(\alpha)$. 

In a factor system $\MMM_\III$ the set $\elem(\MMM_\III)$ is, up to isomorphism, the limit $\limf(\MMM_\III)$, whereby projections have been defined by (\ref{Projdef}) and embeddings are
\begin{align}
\label{Embdef}
Emb_i(a_i) &:= \{\emb{i}{i'}(a_{i}) \in \bigcup_{j \in \III} \MMM_j \mid i' \geq i\}^m.
\end{align} 

Let $\MMM$ be a target for a system $\MMM_\III$, then $\MMM$ is \emph{maximal (over $\MMM_\III$)} iff $Ext(a) \in \elem(\MMM_\III)$ for all $a \in \MMM$. $\MMM$ is \emph{extensional (over $\MMM_\III$)} iff $Ext(a) \sim Ext(b)$ implies $a = b$ for all $a, b \in \MMM$. $\MMM$ is \emph{complete (over $\MMM_\III$)} iff for all consistent sets $\alpha$ there is a limit element $a$, i.e., some $a \in \MMM$ with $Ext(a) \sim \alpha$. One can characterize $\limf(\MMM_\III)$ also as a target that is maximal, extensional and complete over $\MMM_\III$.

\subsection{Stronger Properties for $\fil$}

Condition (\ref{cofineq}) suffices to prove the basic properties of targets and limits of factor systems \cite{eberl2023}. For instance, if $\MMM$ and $\NNN$ are targets for the factor systems $\MMM_\III$ and $\NNN_\JJJ$ respectively, then the set $[\MMM \to_{\fil} \NNN] := \{f : \MMM \to \NNN \mid \III_f \in \fil(\III \times \JJJ)\}$ is a target for the factor system $[\MMM_\III \to \NNN_\JJJ]$. To prove a corresponding statement for limits, two more conditions are needed. The first is a property that relates the available indices of functions to those of their arguments:
\begin{equation}
\label{csetcond}
\tag{D}
\HHH \in \fil(\III \times \JJJ) \text{ and } \III' \in \fil(\III) \text{ implies } \HHH[\III'] \in \fil(\JJJ)
\end{equation}
for all sets $\III$ and $\JJJ$, whereby $\HHH[\III'] := \{j \in \JJJ \mid \exists i \in \III' \ i \to j \in \HHH\}$. This condition guarantees that all functions, given only by their approximations, can be applied to their arguments, given also only by their approximations. Let $\zeta$ be a consistent set on the factor system $[\MMM_\III \to \NNN_\JJJ]$, and $\alpha$ a consistent set on the factor system $\MMM_\III$. The definition
\begin{equation}
\label{appcond}
\tag{Appl}
\zeta(\alpha) := \{f_{i \to j}(a_i) \mid f_{i \to j} \in \zeta \text{ and } a_i \in \alpha\}
\end{equation}
defines an application on consistent sets. We need Condition (\ref{csetcond}) to show that the result is a consistent set again. Moreover, a strengthening of (\ref{cofineq}) is necessary as well:
\begin{equation}
\label{filtereq}
\tag{Filter}
\fil(\III) \text{ is a filter on } \{\HHH \subseteq \fil(\III) \mid \HHH \text{ is cofinal}\}, \text{ containing all up-sets}.
\end{equation}

By filter we always mean a proper filter. Both conditions, (\ref{csetcond}) and (\ref{filtereq}), guarantee that the set $[\MMM \to_{\fil} \NNN]$ is a limit of $[\MMM_\III \to \NNN_\JJJ]$, provided $\MMM$ and $\NNN$ are limits of the factor systems $\MMM_\III$ and $\NNN_\JJJ$ resp. So we have
\begin{equation}
\label{elemisomeq}
[\elem(\MMM_\III) \to_{\fil} \elem(\NNN_\JJJ)] \simeq \elem([\MMM_\III \to \NNN_\JJJ]).
\end{equation}

For the rest of this paper we assume that both conditions, (\ref{csetcond}) and (\ref{filtereq}), are met.

\section{Sets with Families of PERs}
\label{peruniv}

Factor systems $\MMM_\III$ consist of approximations, such as the maps $id_{i \to j} : \nat_i \to \nat_j$, where $n \mapsto n$ for $j \geq i$ are approximations of the identity function $\lambda \y{x}^{nat} \, \y{x} : nat \to nat$. The completed identity function $id : \nat \to \nat$ is not part of the factor system. It is part of the factor limit if one accepts actual infinite sets, including the function space $[\nat \to \nat]$. If limits are seen as indefinitely large finite sets within the system, then the limit set is $[\nat_{i'} \to \nat_{j'}]$ for some indefinitely large index $i' \to j'$.

In both cases, limits $\MMM$ of factor systems $\MMM_\III$ carry a structure that allows us to define a notion of continuity independently of the underlying factor system and independently of the concrete structure of the specific limit $\elem(\MMM_\III)$: Each target $(\MMM, \Pmap)$ for a system $(\MMM_\III, \pmap)$ defines a family of symmetric relations $\per_i$ on $\MMM$ by
\begin{equation}
\label{perpmap}
\tag{Equiv}
a \per_i b \ \iffdef a \Pmap a_i \text{ and } b \Pmap b_i \text{ for some } a_i \comp b_i \text{ in } \MMM_i,
\end{equation}
which satisfy $a \per_i b \imp a \per_i a$ for all $a,b \in \MMM$. If the system is a prefactor system, then each $\per_i$ is transitive as well, i.e., it is a PER (a partial equivalence relation). The idea is that each relation $\per_i$ is an approximation of the equality, and $\per_{i'}$ provides a better approximation than $\per_i$ for $i' \geq i$. The challenge is the partiality of the relations $\per_i$. This partiality is a consequence of the fact that an element does not need to have an approximation at all stages, e.g.~$\lambda \y{x}^{nat} \, \y{x} : nat \to nat$ has only approximations $id_{i \to j} : \nat_i \to \nat_j$ if $j \geq i$.

\subsection{Approximating Equality and the Internal System}
\label{prelimsec}

As in Section \ref{factlimsec}, let $\III$, $\JJJ$ and $\HHH$ be non-empty, directed index sets. Let $\MMM$ be a set endowed with a family $\per_\III \, := (\per_i)_{i \in \III}$ of PERs $\per_i \, \subseteq \MMM \times \MMM$. The set $\{b \in \MMM \mid b \per_i a\}$ is abbreviated by $[a]_i$, and $[i]$ stands for the domain of $\per_i$, i.e., for the set $\{a \in \MMM \mid a \per_i a\}$. Define
\[
\III_a := \{i \in \III \mid a \in [i]\}
\]
for $a \in \MMM$, so $\III_a$ is the set of stages at which $a$ has an approximation. Note that for an element $a$ in a target $\MMM$, the set $\III_a$ from Section \ref{tarsec} conforms to this definition if (\ref{perpmap}) is taken as the basis, which we always do in the following. An immediate consequence of (\ref{perpmap}) is the fact that a target $(\MMM, \Pmap)$ for a system $(\MMM_\III, \pmap)$ satisfies for $a \in \MMM$
\begin{equation}
\label{perpmapcons}
i \in \III_a \iff a \in [i] \iff \exists a_i \in \MMM_i \ a \Pmap a_i \iff \exists a_i \in \MMM_i \ a_i \in Ext(a).
\end{equation}

The next definition introduces the very basic properties that allow us to consider $\per_\III$ as approximations of the equality.

\begin{definition}[$\fil$-set]
\label{perunivdef}
$(\MMM,\per_\III)$ is a $\fil$-set iff 
\begin{enumerate}
\item $\per_i \, \subseteq \MMM \times \MMM$ is a PER for all $i \in \III$,
\item $\per_\III$ is \emph{dense}, i.e., $\III_a \in \fil(\III)$ for all $a \in \MMM$, and
\item $\per_\III$ is \emph{approximating}, i.e., $a \per_{i'} b$ and $a, b \in [i]$ implies $a \per_i b$ for all $i \leq i'$ and all $a, b \in \MMM$.
\end{enumerate}
The family $\per_\III$, and likewise the $\fil$-set $(\MMM,\per_\III)$, is \emph{extensional} iff for all $a, b \in \MMM$, whenever $a \per_i b$ for cofinal many $i \in \III$, then $a = b$. 
\end{definition}

The requirement that $\per_\III$ is dense implies that all elements in $\MMM$ can be identified arbitrarily exact, since all sets in $\fil(\III)$ are cofinal. Note however that two elements $a, b \in \MMM$ need not be equal if $a \per_i b$ holds for cofinal many indices $i \in \III$ --- this is the case only if $\MMM$ is additionally extensional. 

\begin{lemma}[Target structure of a prefactor system]
\label{approxlem}
Let $(\MMM,\Pmap)$ be a target for a prefactor system $(\MMM_\III,\pmap)$, then $(\MMM,\per_\III)$, with $\per_\III$ as defined in (\ref{perpmap}), is a $\fil$-set.
\end{lemma}

\begin{proof}
Density follows immediately from (\ref{perpmapcons}). To see that $\per_\III$ is approximating, let $a \per_{i'} b$. Then there are $a_{i'}, b_{i'} \in \MMM_{i'}$ with $a \Pmap a_{i'}$, $b \Pmap b_{i'}$ and $a_{i'} \comp b_{i'}$. For $a, b \in [i]$ there are moreover $a_{i}, b_{i} \in \MMM_{i}$ with $a \Pmap a_{i}$ and $b \Pmap b_{i}$. We claim that $a_i \comp b_i$. From (\ref{ppmapeq}) we infer $a_{i'} \pmap a_i$ and $b_{i'} \pmap b_i$. Now $a_{i'} \comp b_{i'} \pmap b_i$ implies $a_{i'} \pmap b_i$, hence $a_i \comp b_i$ by definition.
\end{proof}

\begin{example}
\label{conssetex}
The set of all consistent sets of a prefactor system $\MMM_\III$ form a $\fil$-set. Definition (\ref{perpmap}) gives $\alpha \per_i \beta \iff a_i \comp b_i$ for some $a_i \in \alpha$ and $b_i \in \beta$. The same holds for the set of maximal consistent sets, i.e., for dynamic elements $\elem(\MMM_\III)$. 
\end{example}

The next definition of a homomorphism already contains the enriched structures that we introduce later.

\begin{definition}
\label{homdef2}
A \emph{homomorphism} $\Phi = (\Phi_0, \Phi_1)$ between two $\fil$-sets $(\MMM, \per_\III)$ and $(\NNN, \per_\JJJ)$ consists of maps $\Phi_0 : \III \to \JJJ$ and $\Phi_1 : \MMM \to \NNN$ such that $\Phi_0$ is monotone and for all $i \in \III$
\begin{equation}
\label{homeq2}
a \per_i b \ \imp \ \Phi_1(a) \per_j \Phi_1(b)
\end{equation}
for $j := \Phi_0(i)$. If the equivalence holds in (\ref{homeq2}), then $\Phi$ is said to be \emph{strong}. A homomorphism between two pointed $\fil$-sets $(\MMM, \per_\III, pt_\III)$ and $(\NNN, \per_\JJJ, pt_\JJJ)$, introduced in Definition \ref{pointperdef}, additionally satisfies $\Phi_1(pt_{i}(a)) = pt_j(\Phi_1(a))$. A homomorphism between PER-sets $(\MMM, \per_\III, \ptp_\III)$, introduced in Definition \ref{extperdef}, satisfies the stronger property $\Phi_1(\ptp_{i}(a)) = \ptp_j(\Phi_1(a))$.

We call $\Phi$ injective (surjective) iff both maps $\Phi_0$ and $\Phi_1$ are injective (surjective resp.). An \emph{isomorphism} $\Phi$ between two $\fil$-sets is a homomorphism with a further such homomorphism $\Psi$ as its inverse, i.e., $\Psi_0$ and $\Psi_1$ are inverse to $\Phi_0$ and $\Phi_1$ resp. So an isomorphism is automatically strong. We write $\MMM \simeq \NNN$ if an isomorphism between $(\MMM, \per_\III)$ and $(\NNN, \per_\JJJ)$ exists.
\end{definition}

If $\III = \JJJ$ and $\Phi_0$ is the identity map, we do not mention $\Phi_0$. This will be the situation in the sequel.

\begin{example}
\label{Exthomoex}
Let $\MMM$ be a target for a prefactor system $\MMM_\III$. The map $Ext : \MMM \ni a \mapsto Ext(a)$ is a strong homomorphism between the $\fil$-sets $\MMM$ (see Lemma \ref{approxlem}) and the set of all consistent sets on $\MMM_\III$ (see Example \ref{conssetex}). Similarly,
\[
Ext^m : \MMM \to \elem(\MMM_\III), \ a \mapsto Ext(a)^m
\]
is a homomorphism which is strong if $\MMM$ is maximal over $\MMM_\III$.
\end{example}

In the subsequent sections we will always assume that at each stage $i$ there are only finitely many equivalence classes.\footnote{This assumption is used only in order to show that the related factor systems have finitely many objects at each stage, all other concepts do not use it. If one is not interested in a finitistic view, this assumption can be dropped here and for systems as well, or replaced by a more general notion of ``definite'' sets.} 
 
\begin{definition}[Internal system]
\label{intsysdef}
Define $\APX(\MMM) := ((\quot{\MMM}{\per_{i}})_{i \in \III},\pmap)$ with $\pmap$ given by overlapping equivalence classes, i.e.,
\begin{equation}
\label{intsyseq}
A_{i'} \pmap A_i \ \iffdef A_{i'} \cap A_i \not= \emptyset
\end{equation}
for $A_{i'} \in \quot{\MMM}{\per_{i'}}$, $A_i \in \quot{\MMM}{\per_{i}}$ and $i' \geq i$. $\APX(\MMM)$ is called the \emph{internal system} of $(\MMM,\per_\III)$.
\end{definition}

$\APX(\MMM)$ is indeed a system since relation $\pmap$ is obviously reflexive and there are only finitely many equivalence classes at each stage $i \in \III$ by assumption.

\begin{lemma}[Internal system of a $\fil$-set]
\label{densefamlem}
The internal system $\APX(\MMM)$ of a $\fil$-set $(\MMM,\per_\III)$ is a system which satisfies (\ref{idmapeq}). $(\MMM,\in)$ is a target for the internal system, that is, $\Pmap$ is defined as\footnote{Compare this with (\ref{elemtergeteq}), where the membership relation has the opposite direction when the limit is constructed from the system.}
\[
a \Pmap A_i \, \iffdef a \in A_i.
\]

Then $Ext(a) = \{A_i \in \bigcup_{i \in \III} \quot{\MMM}{\per_{i}} \mid a \in A_i\}$ for $a \in \MMM$ and $a \per_i b \iff a, b \in A_i$ for some $A_i \in \quot{\MMM}{\per_{i}}$, so the definition of $\Pmap$ is in line with (\ref{perpmap}).
\end{lemma}

\begin{proof}
Property (\ref{idmapeq}) holds because $A_i \comp B_i \iff A_i \cap B_i \not= \emptyset \iff A_i = B_i$. Being a target of a system requires that $\{i \in \III \mid \exists A_i \in \MMM_i \ a \Pmap A_i\} = \{i \in \III \mid a \in [i]\} \in \fil(\III)$ for all objects $a \in \MMM$, which holds by density. The additional claims are easily verified.
\end{proof}

Maximality, extensionality and completeness have been defined for targets for prefactor systems, see Section \ref{tarsec}. In Definition \ref{perunivdef} we introduced a notion of extensionality for a $\fil$-set $(\MMM,\per_\III)$ which does not refer to a system, but uses only the PERs $\per_i$. It is easy to see, and proven in Lemma \ref{extcorlem}, that the $\fil$-set $(\MMM,\per_\III)$ is extensional iff $(\MMM,\in)$ is extensional over the prefactor system $\APX(\MMM)$. As far as the other properties are concerned, $(\MMM,\Pmap)$ is automatically maximal over $\APX(\MMM)$ when we introduce points in Section \ref{richuni}, and completeness is the subject of Section \ref{complpersec}. 
\begin{lemma}
\label{extcorlem}
Given a target $(\MMM, \Pmap)$ for a prefactor system $(\MMM_\III,\pmap)$. Then the following are equivalent:
\begin{enumerate}
\item \label{extcorlem1} $(\MMM, \Pmap)$ is extensional over $\MMM_\III$,
\item \label{extcorlem2} $(\MMM, \per_\III)$ is extensional,
\item \label{extcorlem3} for all $a, b \in \MMM$, whenever $a \per_i b$ for all $i \in \III_a \cap \III_b$, then $a = b$,
\item \label{extcorlem4} $(\MMM,\in)$ is extensional over $\APX(\MMM)$.
\end{enumerate}
\end{lemma}

\begin{proof}
The equivalence of \ref{extcorlem1}.~and \ref{extcorlem3}.~is established if we prove that $Ext(a) \sim Ext(b)$ holds iff $a \per_i b$ for all $i\in \III$ and all $a, b \in [i]$. Now $Ext(a) \sim Ext(b)$ is equivalent to $a_i \comp b_i$ for all $a_i \in Ext(a)$ and $b_i \in Ext(b)$, which again is equivalent to $a \per_i b$ by (\ref{perpmap}).

In order to establish the equivalence between \ref{extcorlem2}.~and \ref{extcorlem3}.~it suffices to show that if $a \per_i b$ for cofinal many $i \in \III_a \cap \III_b$, then $a \per_i b$ for all $i \in \III_a \cap \III_b$. So given any $i \in \III_a \cap \III_b$. By assumption there is some $i' \geq i$, $i' \in \III_a \cap \III_b$ such that $a \per_{i'} b$. Since $\per_\III$ is approximating by Lemma \ref{approxlem}, we conclude that $a \per_{i} b$.

The equivalence of \ref{extcorlem3}.~and \ref{extcorlem4}.~follows from the fact that $\per_i$ is indeed a PER for a prefactor system $\MMM_\III$. So for the internal system we have $Ext(a) \sim Ext(b)$ iff $\{A_i \in \bigcup_{i \in \III} \quot{\MMM}{\per_{i}} \mid a \in A_i\} \sim \{B_i \in \bigcup_{i \in \III} \quot{\MMM}{\per_{i}} \mid b \in B_i\}$. This is the case iff $A_i = B_i$ for all $A_i \in Ext(a)$, $B_i \in Ext(b)$ and $A_i, B_i \in \quot{\MMM}{\per_{i}}$ --- note that $A_i \comp B_i \iff A_i = B_i$ by Lemma \ref{densefamlem}. This translates to $a, b \in A_i \in \quot{\MMM}{\per_{i}}$, which means $a \per_i b$.
\end{proof}

The same reasoning shows the equivalences in Lemma \ref{extcorlem} also for a prefactor system with embeddings/projections as well as for a factor system.

\subsection{Pointed Families of PERs}
\label{richuni}

It is possible to define a property for $\per_\III$ that corresponds to the basic Property (\ref{pmapeq}) of a prefactor system, but this property has no intuitive reading. For that reason, and for the sake of simplicity, we will directly consider \emph{prefactor systems with embeddings}. The counterpart to the embeddings in a factor system are points (specific elements) in each of the equivalence classes of the PER. So we introduce pointed families of PERs. A point is, roughly, a canonical element in the equivalence class with the property that if the precision increases, then the set of these canonical elements increases as well, so they do not ``switch'' during the expansion.

\begin{definition}[Pointed $\fil$-set]
\label{pointperdef}
A family of PERs $\per_\III$ on $\MMM$ is \emph{pointed} (or \emph{has points}) iff there are sets $\PPP_i \subseteq \MMM$ for all $i \in \III$ and surjective maps $pt_i : [i] \to \PPP_i$ such that for all $a, b \in \MMM$ and all $i \leq i'$
\begin{enumerate}
\item \label{2pointperdef} $\PPP_i \subseteq \PPP_{i'}$,
\item \label{1pointperdef} $a \per_i b$ implies $a \per_i pt_i(a) = pt_i(b)$, 
\item \label{3pointperdef} $a \in [i]$ implies $pt_{i'}(a) \in [i]$ provided that $a \in [i']$.
\end{enumerate}
Sometimes we call the elements in $\PPP_i$ \emph{grid points}. A $\fil$-set $(\MMM,\per_\III)$ is pointed iff $\per_\III$ is pointed.
\end{definition}

A pointed $\fil$-set is thus given by $(\MMM,\per_\III,pt_\III)$ with $pt_\III = (pt_i)_{i \in \III}$. One may think of a pointed family of PERs as providing an identification of objects in terms of their currently known value (to the precision of $i$) as well as possible more precise values that may be known in the future. Elements in $\PPP_i$ describe the current values of objects given with precision $i$, i.e., their states at $i$. The equivalence class $[a]_i$ contains all possible values that we could discover for an object $a$, if we identify it with a precision $i' \geq i$. Next we state some simple, useful facts:

\begin{lemma}
\label{simplepointlem}
Given a pointed family of PERs $(\per_\III,pt_\III)$ on $\MMM$ and indices $i \leq i'$, then
\begin{enumerate}
\item $\PPP_i \subseteq [i']$,
\item $a \per_i b$ iff $pt_i(a) = pt_i(b)$ for all elements $a, b \in [i]$,
\item $pt_i(b) = a$ iff $a \in \PPP_i$ and $a \per_i b$,
\item $pt_{i'}(a) = a$ for all $a \in \PPP_i$,
\item $pt_{i'}(pt_{i}(a)) = pt_{i}(a)$ for all $a \in [i]$.
\end{enumerate}
\end{lemma}

\begin{proof}
First we shall prove that for all $a \in \PPP_i$ we have $a \in [i]$ and $pt_{i}(a) = a$: If $a \in \PPP_i$, then there is an element $b \in [i]$ with $a = pt_i(b)$. Because $pt_i(b) \per_i b$, that is $a \per_i b$, we conclude that $a \in [i]$ as well as $pt_i(a) = pt_i(b) = a$.

The first relation follows thus from $\PPP_i \subseteq \PPP_{i'} \subseteq [i']$. For the backward implication of the second equivalence note that $a \per_i pt_i(a) = pt_i(b) \per_i b$ holds for $a, b \in [i]$, thus $a \per_i b$. For the third part assume $pt_i(b) = a$, then clearly $a \in \PPP_i$ and $a = pt_i(b) \per_i b$. Conversely we have $pt_i(b) = pt_i(a) = a$, the latter equation holds since $a \in \PPP_i$. For the forth part note that $a \in \PPP_i \subseteq \PPP_{i'}$ implies $pt_{i'}(a) = a$. This also establishes the last part --- simply observe that $pt_i(a) \in \PPP_i$ for $a \in [i]$. 
\end{proof}

If in Definition \ref{pointperdef} the last part is replaced by the stronger property $pt_{i'}(a) \per_i a$, or equivalently,
\begin{equation}
\label{approxpointeq}
pt_i(pt_{i'}(a)) = pt_i(a) \text{ for all } a \in [i] \cap [i'] \text{ and } i \leq i',
\end{equation}
then the family of PERs is automatically approximating. More precisely, we have:

\begin{lemma}
\label{eqivpointlem}
Given a family of PERs $\per_\III$ with maps $pt_i : [i] \to \PPP_i$ satisfying the first two parts of Definition \ref{pointperdef}. Then part \ref{3pointperdef} and being approximating together are equivalent to (\ref{approxpointeq}). 
\end{lemma}

\begin{proof}
Assume that $\per_\III$ is approximating and part \ref{3pointperdef} holds. If $a \in [i] \cap [i']$, then $a, pt_{i'}(a) \in [i]$. Hence $pt_{i'}(a) \per_i a$ follows from $pt_{i'}(a) \per_{i'} a$, since $\per_\III$ is approximating. 

Conversely, assume (\ref{approxpointeq}), which clearly implies part \ref{3pointperdef}, so let $a \per_{i'} b$ with $a, b \in [i]$ be given. By assumption $a \per_i pt_{i'}(a)$ and $b \per_i pt_{i'}(b)$, thus $a \per_i pt_{i'}(a) = pt_{i'}(b) \per_i b$, which shows that $\per_\III$ is approximating.
\end{proof}

\begin{lemma}
\label{eqivpointlem1}
For a pointed $\fil$-set we have 
\[
pt_{i'}(pt_{i}(a)) = pt_{i}(pt_{i'}(a)) = pt_{i}(a)
\]
for $a \in [i] \cap [i']$ and $i \leq i'$.
\end{lemma}

\begin{proof}
We claim that $pt_{i}(pt_{i'}(a)) = pt_{i}(a)$ --- the equation $pt_{i'}(pt_{i}(a)) = pt_{i}(a)$ has already be shown in Lemma \ref{simplepointlem}. Indeed, $pt_{i'}(a) \per_{i'} a$ implies $pt_{i'}(a) \per_{i} a$ since $\per_\III$ is approximating and $a, pt_{i'}(a) \in [i]$. This establishes $pt_{i}(pt_{i'}(a)) = pt_{i}(a)$ as claimed.
\end{proof}

Each target for a prefactor system with embeddings generates in a natural way a pointed $\fil$-set --- the proof is postponed to the appendix:

\begin{proposition}[Target structure of a prefactor system with embeddings]
\label{targpointlem}
Assume we have a target $(\MMM,\Pmap,Emb)$ for a prefactor system with embeddings $(\MMM_\III, \pmap, emb)$. Then $(\MMM,\per_\III,pt_\III)$, where $\per_i$ is defined by (\ref{perpmap}), $\PPP_i := \{Emb_i(a_i) \mid a_i \in \MMM_i\}$ and
\begin{equation}
\label{perpmap2}
\tag{Point}
pt_i(a) = b \, \iffdef a \Pmap a_i \text{ and } b = Emb_i(a_i) \text{ for some } a_i \in \MMM_i
\end{equation}
is a pointed $\fil$-set. In particular,
\begin{enumerate}
\item \label{2targpointlem} $a_{i} \comp b_i \iff Emb_{i}(a_{i}) = Emb_i(b_i)$ and
\item \label{3targpointlem} $a_{i'} \pmap a_i \, \imp \, pt_i(Emb_{i'}(a_{i'})) = Emb_i(a_i)$, and if $\MMM_\III$ is stable, the inverse implication holds as well.
\end{enumerate}

If $\MMM_\III$ is stable, then $Emb : \MMM_\III \to \APX(\MMM)$ is a strong and surjective homomorphism, and if $\pmap$ satisfies (\ref{idmapeq}), then $Emb$ is an isomorphism.
\end{proposition}

If we speak of a target as a set with a pointed family of PERs, we refer to the implicit structure given by Equations (\ref{perpmap}) and (\ref{perpmap2}).

\subsection{The Internal System for Pointed Families of PERs}
\label{intfactsys}

The internal system $\APX(\MMM)$ has been introduced in Definition \ref{intsysdef}. If we consider pointed families of PERs, we can take the points $\PPP_\III$ instead of the equivalence classes. The property $A_{i'} \cap A_i \not= \emptyset$ from (\ref{intsyseq}) becomes $a_{i'} \per_i a_i$ for the unique points $a_{i'} \in A_{i'} \cap \PPP_{i'}$ and  $a_i \in A_i \cap \PPP_{i}$, so we define for $a_{i'} \in \PPP_{i'}$, $a_i \in \PPP_i$ and $i' \geq i$
\begin{equation}
\label{pmapmultuni}
a_{i'} \pmap a_i \iffdef a_{i'} \per_i a_i \iff pt_i(a_{i'}) = a_i.
\end{equation}

Then $\APX(\MMM) = (\PPP_\III,\pmap)$ has automatically embeddings, being subset inclusion $\PPP_i \embb \PPP_{i'}$. Condition (\ref{embcond}) is thus trivially satisfied. 

\begin{proposition}[Internal system of a pointed $\fil$-set]
\label{apxprelimlem}
Given a pointed $\fil$-set $(\MMM,\per_\III,pt_\III)$. Then the internal system $\APX(\MMM)$ is a prefactor system with embeddings. $\MMM$ together with
\[
a \Pmap a_i \iffdef a_i = pt_i(a)
\]
and $Emb_i : \PPP_i \embb \MMM$, $Emb_i(a) := a$ is a target for $\APX(\MMM)$; for such a target we write $(\MMM,pt_\III,\embb)$. The extension of $a \in \MMM$ is $Ext(a) = (pt_i(a))_{i \in \III_a}$ and $(\MMM,pt_\III,\embb)$ is maximal over $\APX(\MMM)$.
\end{proposition}

\begin{proof}
To show Property (\ref{pmapeq}) assume that $a_{i''} \in \PPP_{i''}$, $a_{i'} \in \PPP_{i'}$, $a_i \in \PPP_{i}$ with $pt_{i'}(a_{i''}) = a_{i'}$ and $pt_i(a_{i''}) = a_i$. Then $pt_i(a_{i'}) = pt_i(pt_{i'}(a_{i''})) = pt_i(a_{i''}) = a_i$ by Lemma \ref{eqivpointlem1}. $\APX(\MMM)$ trivially has embeddings, as noted above.

$(\MMM,pt_\III,\embb)$ is a target for the prefactor system with embeddings $\APX(\MMM)$: We have to show that $a \Pmap a_{i'}$ and $a \Pmap a_{i}$ implies $a_{i'} \pmap a_i$, which however follows from $pt_i(a_{i'}) = pt_i(pt_{i'}(a)) = pt_i(a) = a_i$ by Lemma \ref{eqivpointlem1}. The required Property (\ref{embcond}), i.e., $a' \pmap a \imp Emb_{i'}(a') \Pmap a$ holds trivially.

$\MMM$ is maximal over $\MMM_\III$, that is, $Ext(a) = (pt_i(a))_{i \in \III_a}$ is maximal as a consistent set: Assume $pt_{i'}(a) \pmap b \in \PPP_i$, then by definition $pt_i(pt_{i'}(a)) = b$, consequently $pt_i(a) = b$, which shows $b \in Ext(a)$.
\end{proof}

Proposition \ref{targpointlem} states that a target for a prefactor system with embeddings automatically has the structure of a pointed $\fil$-set. Proposition \ref{apxprelimlem} states the converse in the sense that if we start from a pointed $\fil$-set, then it automatically contains a prefactor system with embeddings such that $\MMM$ naturally gives rise to a target for it. This system and the target are of a certain kind: 
\begin{enumerate}
\item The system satisfies Property (\ref{idmapeq}) by Lemma \ref{densefamlem} and 
\item the target is maximal over the system by Proposition \ref{apxprelimlem}.
\end{enumerate}

The internal system $\APX(\MMM)$ of a target $\MMM$ for $\MMM_\III$ is generally not isomorphic to $\MMM_\III$, because if the Property (\ref{idmapeq}) does not hold, the elements are identified in the step from $\MMM_\III$ to $\MMM$ and cannot be recovered. But if we take ``$\MMM_\III$ modulo $\comp$'', there is a one-to-one correspondence. This is basically the statement about homomorphisms in Proposition \ref{targpointlem}.

\subsection{Convergence and Completeness}
\label{complpersec}

A consistent set $\alpha$ in a system $\MMM_\III$ is basically a convergent sequence. In particular, if Property (\ref{idmapeq}) is satisfied, then $\alpha$ can be organized as a family $(a_i)_{i \in \JJJ}$ for an index set $\JJJ \in \fil(\III)$. Remind that a limit element $a$ of $\alpha$ in a target $\MMM$ of $\MMM_\III$, if it exists, satisfies $\alpha \sim Ext(a)$. If $\MMM = \elem(\MMM_\III)$ is the limit of $\MMM_\III$, then $\alpha^m = Ext(a)$. Now we also define convergent sequences in $\fil$-sets.

\begin{definition}[Convergent sequence of points]
Given a pointed $\fil$-set $(\MMM,\per_\III,pt_\III)$. Then $(a^i)_{i \in \JJJ}$ with $a^i \in \PPP_i$ and $\JJJ \in \fil(\III)$ is a \emph{convergent sequence of points}\footnote{One may define general sequences in $\MMM$ as well, not only sequences of points. We will not do this here, as the partiality of the sequences requires some care. Instead we will concentrate on the construction of the function space in Section \ref{funPERset}.}  iff $pt_i(a^{i'}) = a^i$ for all $i' \geq i$. An element $a \in \MMM$ is called a \emph{limit element} of a convergent sequence of points $(a^i)_{i \in \JJJ}$ iff $pt_i(a) = a^i$ for all $i \in \JJJ$. We denote $a$ by $lim_{i \in \JJJ}a^i$ if this element is unique.
\end{definition}

\begin{figure}[ht]
\centering
\includegraphics[trim = 0mm 0mm 0mm 0mm, width=1\textwidth]{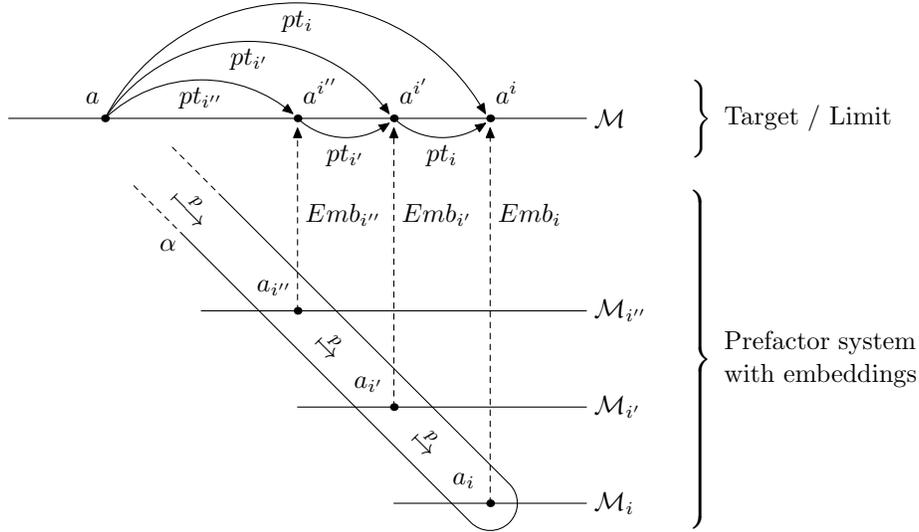}
\caption{Convergent sequences in systems and limits.}
\label{fig1}
\end{figure}

It is easy to see that in the above definition $(a^i)_{i \in \JJJ}$ is a convergent sequence of points iff $(a^i)_{i \in \JJJ}$ is a consistent set in the internal system of $\MMM$. Lemma \ref{comppreplem} formulates a more general statement, outlined in Figure \ref{fig1}:

\begin{lemma}
\label{comppreplem}
Given a target $(\MMM, \Pmap, Emb)$ for a prefactor system with embeddings $(\MMM_\III,\pmap,emb)$. If $\alpha$ is a consistent set in $\MMM_\III$, then $(a^i)_{i \in \III_\alpha}$ with $a^i := Emb_i(a_i) \in \MMM$ for $a_i \in \alpha$ is a convergent sequence of points. If $\MMM_\III$ is stable, the converse holds as well. Moreover, $\alpha \sim Ext(a) \iff pt_i(a) = a^i$ for all $i \in \III_\alpha \cap \III_a$. 
\end{lemma}

\begin{proof}
In the view of part \ref{2targpointlem} of Proposition \ref{targpointlem}, $(a^i)_{i \in \III_\alpha}$ is well defined, since $a_i, b_i \in \alpha$ implies $a_i \comp b_i$, which is then a convergent sequence of points by part \ref{3targpointlem}. Moreover, the statement $\alpha \sim Ext(a)$ is equivalent to $b_i \comp a_i$ for all $b_i \in \alpha$, $a_i \in Ext(a)$ and $i \in \III_a \cap \III_\alpha$. Now $b_i \comp a_i$ is equivalent to $Emb_i(b_i) = Emb_i(a_i)$, hence $a^i = Emb_i(b_i) = Emb_i(a_i) = pt_i(a)$, whereby the last equation holds due to (\ref{perpmap2}).
\end{proof}

The last claim in Lemma \ref{comppreplem} states that $a \in \MMM$ is a limit element of the consistent set $\alpha$ iff it is a limit element of the convergent sequence of points $(a^i)_{i \in \III_\alpha}$, given by $Emb[\alpha]$. And if the limit element is a unique, then $lim(\alpha) = lim_{i \in \III_\alpha}a^i$. For this reason, similar to extensionality, we can prove equivalent conditions for completeness. Note that the next lemma also holds for targets for stable factor systems. 

\begin{lemma}
\label{limcorlem}
Given a target $(\MMM, \Pmap, Emb)$ for a stable prefactor system with embeddings $(\MMM_\III,\pmap,emb)$. Then the following are equivalent:
\begin{enumerate}
\item $(\MMM, \Pmap, Emb)$ is complete over $\MMM_\III$,
\item every convergent sequence of points has a limit element,
\item $(\MMM, \in, \embb)$ is complete over $\APX(\MMM)$.
\end{enumerate}
\end{lemma}

\begin{proof}
By Lemma \ref{comppreplem} the consistent sets $\alpha$ in $\MMM_\III$ correspond one-to-one to the convergent sequences of points $Emb[\alpha]$. For $\alpha$ exists an element $a \in \MMM$ with $\alpha \sim Ext(a)$ exactly if $a$ is limit of this convergent sequence, again by Lemma \ref{comppreplem}. This establishes the equivalence between the first to claims. The last statement is equivalent to part 2, since $(\MMM, \in, \embb)$ is a target for its internal system and we may apply the same argument as before.
\end{proof}

This leads us to the following definition:

\begin{definition}
\label{extcompldef}
A pointed $\fil$-set $(\MMM,\per_\III,pt_\III)$ is \emph{complete} iff every convergent sequence of points has a limit element.
\end{definition}

Thus a pointed $\fil$-set $(\MMM,\per_\III,pt_\III)$ is complete iff it is complete over a stable prefactor system with embeddings $(\MMM_\III,\pmap,emb)$. A similar statement for extensionality was proven in Lemma \ref{extcorlem}. So we can use the terms \emph{extensional} and \emph{complete} in this situation, without mentioning whether they apply to $\MMM$ as a target or as a pointed $\fil$-set.

A prefactor system with embeddings $\MMM_\III$ defines a simultaneous process of sets and their elements, where the construction of $\elem(\MMM_\III)$ from $\MMM_\III$ provides a notion of a limit for both, sets and elements. Thus they are a candidate for a potential infinite model together with an interpretation of types and terms in it. The only shortcoming so far is that they are not closed under the function space construction.

\subsection{Total Equivalence Relations and PER-Sets}
\label{eqrelpersec}

For the function space construction of sets with families of PERs the counterpart of projections in a factor system is required, i.e., the PERs $\per_i$ are extended to (total) equivalence relations $\perp_i$. In the case of PERs with points, the extended equivalence relations $\perp_i$ are a consequence of an extension of $pt_i$ to total functions $\ptp_i$. We will now define the extension of $pt_\III$ to $\ptp_\III$, which shall preserve basic properties of $pt_\III$.

\begin{definition}[PER-set]
\label{extperdef}
Given a pointed $\fil$-set $(\MMM,\per_\III,pt_\III)$. Then $\ptp_\III$ \emph{extends} $pt_\III$ iff $\ptp_i : \MMM \to \PPP_i$ and for all $i' \geq i$
\begin{enumerate}
\item $\ptp_i(a) = pt_i(a)$ for $a \in [i]$,
\item $\ptp_i = \ptp_i \circ \ptp_{i'}$ and
\item $a \in [i]$ implies $\ptp_{i'}(a) \in [i]$.
\end{enumerate}

Define $a \perp_{i} b$ as $\ptp_i(a) = \ptp_i(b)$ for $a, b \in \MMM$. A \emph{PER-set} is a pointed $\fil$-set with a function $\ptp_\III$ which extends $pt_\III$.\footnote{A PER-set should not be confused with the setoid model of type theory, which has PERs as well, see e.g.~\cite{hofmann1995simple}.}
\end{definition}

A PER-set is given as a structure $(\MMM,\per_\III,\ptp_\III)$. Since $pt_i$ is the restriction of $\ptp_i$ to $[i]$, it is not necessary to mention it. 

\begin{example}
\label{natpointex}
Recall the standard model of the natural numbers from Example \ref{natex} and define $\fil(\nat^+) := \{\HHH \subseteq \nat^+ \mid \up i \subseteq \HHH \text{ for some } i \in \nat^+\}$ with $\up i := \{i' \in \nat^+ \mid i' \geq i\}$. The PER-set structure with points and extension on the limit $\nat$ is:
\begin{enumerate}
\item $n \per_i m \iff n = m \text{ and } n,m < i \in \nat^+$,
\item $\PPP_i = \nat_i$,
\item $Pt_i(n) = \min(n,i-1)$ and thus $n \perp_i m \iff n = m \text{ or both } n,m \geq i-1$.
\end{enumerate}
\end{example}

Let us summarize the properties of a PER-set in a condensed form.

\begin{proposition}
\label{persetdefprop}
Given a set $\MMM$ with a family of PERs $\per_\III \, \subseteq \MMM \times \MMM$ and surjective maps $\ptp_i : \MMM \to \PPP_i \subseteq \MMM$. Let $pt_i$ be the restriction of $Pt_i$ to $[i]$. Then $(\MMM,\per_\III,\ptp_\III)$ is a PER-set iff for all $i \leq i'$
\begin{enumerate}
\item \label{1persetdefprop} $\per_\III$ is dense,
\item \label{2persetdefprop} $\PPP_i \subseteq \PPP_{i'}$,
\item \label{3persetdefprop} $a \per_i b$ implies $a \per_i \ptp_i(a) = \ptp_i(b)$,
\item \label{4persetdefprop} $\ptp_{i} = \ptp_{i} \circ \ptp_{i'}$ and
\item \label{5persetdefprop} $a \in [i]$ implies $\ptp_{i'}(a) \in [i]$.
\end{enumerate}
\end{proposition}

\begin{proof}
If $(\MMM,\per_\III,\ptp_\III)$ is a PER-set, then the parts \ref{1persetdefprop}.~to \ref{5persetdefprop}.~hold by definition, so assume these five properties. We start with Definition \ref{perunivdef} and it suffices to show that $\per_\III$ is approximating. Given $a,b \in [i]$ with $a \per_{i'} b$. We infer with the third and forth part as follows: 
\[
a \per_{i'} b \imp Pt_{i'}(a) = Pt_{i'}(b) \imp Pt_{i}(Pt_{i'}(a)) = Pt_{i}(Pt_{i'}(b)) \imp Pt_{i}(a) = Pt_{i}(b).
\]

Now $a,b \in [i]$ implies $a \per_i a$ and $b \per_i b$, and again with part \ref{3persetdefprop} we conclude that $a \per_i \ptp_{i}(a) = \ptp_{i}(b) \per_i b$, showing that $\per_\III$ is approximating. It is immediate to see that $(\MMM,\per_\III,\ptp_\III)$ is pointed and that the properties of Definition \ref{extperdef} hold.
\end{proof}

The next lemma lists some simple, basic facts.

\begin{lemma}
\label{extperlem}
Given a PER-set $(\MMM,\per_\III,Pt_\III)$, then for all $i \leq i'$,
\begin{enumerate}
\item \label{1extperlem} $a \perp_{i'} b$ implies $a \perp_{i} b$,
\item \label{2extperlem} $a \per_i b$ iff $a \perp_i b$ and $a, b \in [i]$,
\item \label{3extperlem} $a \perp_i \ptp_{i'}(a)$, and $a \per_i \ptp_{i'}(a)$ if $a \in [i]$,
\item \label{5extperlem} $a \per_i b$ implies $\ptp_{i'}(a) \per_i \ptp_{i'}(b)$,
\item \label{6extperlem} $\ptp_{i'}(a) = a$ if $a \in \PPP_i$,
\item \label{7extperlem} $\ptp_{i'} \circ \ptp_{i} = \ptp_{i} \circ \ptp_{i'} = \ptp_{i}$,
\item \label{8extperlem} $pt_{i}(\ptp_{i'}(a)) = pt_{i}(a)$ if $a \in [i]$.
\end{enumerate}
\end{lemma}

\begin{proof}
Part \ref{1extperlem} follows from $\ptp_{i} = \ptp_{i} \circ \ptp_{i'}$ and part \ref{2extperlem} is obvious. For part \ref{3extperlem} it suffices to observe that $a \perp_{i'} \ptp_{i'}(a)$ implies $a \perp_i \ptp_{i'}(a)$ by part \ref{1extperlem}, and that $\ptp_{i'}(a) \in [i]$ if $a \in [i]$. This immediately implies part \ref{5extperlem}, since $\ptp_{i'}(a) \per_i a \per_i b \per_i \ptp_{i'}(b)$. Part \ref{6extperlem} follows from Lemma \ref{simplepointlem} and for part \ref{7extperlem} it suffices to show $\ptp_{i'}(\ptp_{i}(a)) = \ptp_{i}(a)$; this is however a consequence of $\ptp_{i}(a) \in \PPP_i \subseteq \PPP_{i'}$, and the last part is readily seen.
\end{proof}

The PER-sets give rise to enhancements of the propositions already stated for pointed $\fil$-sets. We start with Proposition \ref{targpointlem} --- note that if $\MMM$ is a target for a factor system $\MMM_\III$, it is also a target for $\MMM_\III$ as a prefactor system with embeddings.

\begin{proposition}[Target structure of a factor system]
\label{projextlem}
Let a target $(\MMM,\Pmap,Emb,Proj)$ for a factor system $(\MMM_\III, \pmap, emb, proj)$ be given. Then $(\MMM,\per_\III,\ptp_\III)$, with $\per_i$ as defined in (\ref{perpmap}) together with 
\[
\ptp_i := Emb_i \circ Proj_i
\]
is a PER-set, which satisfies 
\begin{equation}
\label{projembeq}
Proj_i(Emb_i(a_i)) \comp a_i \text{ for all } a_i \in \MMM_i.
\end{equation}

If $\MMM_\III$ is stable, then $Emb : \MMM_\III \to \APX(\MMM)$ is a strong and surjective homomorphism.
\end{proposition}

When we speak of a PER-set as a target for the factor system, we are referring to this automatically included structure of $\per_\III$ and $\ptp_\III$. In the following enhancement of Proposition \ref{apxprelimlem} the notion $\ptp_i \! \restr_{\PPP_{i'}}$ refers to the restriction of $\ptp_i$ to $\PPP_{i'}$:

\begin{proposition}[Internal system of a PER-set]
\label{apxpersetlem}
Given a PER-set $(\MMM,\per_\III,\ptp_\III)$. Then the internal system $\APX(\MMM)$ is a factor system, which has $\proj{i'}{i} := \ptp_i \! \restr_{\PPP_{i'}}$ as projections. $(\MMM,pt_\III,\embb,\ptp_\III)$ is a target for this factor system.
\end{proposition}

\begin{proof}
It is immediate that the maps $\proj{i'}{i}$ are $\comp$-projections and inverse to the embeddings. To prove Property (\ref{projcond}) means showing that $a_{i''}\per_i a_i$ implies $\ptp_{i'}(a_{i''}) \per_i a_i$ for all $i \leq i' \leq i''$. We deduce from part \ref{3extperlem} in Lemma \ref{extperlem} that $\ptp_{i'}(a_{i''}) \per_i a_{i''}$ since $a_{i''} \in [i]$, hence $\ptp_{i'}(a_{i''}) \per_i a_{i''} \per_i a_i$. To see that $(\MMM,pt_\III,\embb,\ptp_\III)$ is a target for the factor system $\APX(\MMM)$ it suffices to check the additional properties of $\ptp_\III$. They follow immediately from parts \ref{6extperlem}, \ref{7extperlem} and \ref{8extperlem} of Lemma \ref{extperlem}.
\end{proof}

The next corollary summarizes a basic connection between standard factor systems and PER-sets:

\begin{corollary}
\label{apxlem}
Given a factor system $\MMM_\III$ which satisfies the Property (\ref{idmapeq}) and let the PER-set $\MMM$ be a maximal target for $\MMM_\III$.
\begin{enumerate}
\item Then $\MMM_\III \simeq \APX(\MMM)$, in particular we have 
\[
\MMM_\III \simeq \APX(\elem(\MMM_\III)).
\]
\item If $\MMM$ is extensional and complete, then $\MMM$ is a limit of $\MMM_\III$ and $\MMM \simeq \elem(\MMM_\III)$, in particular we have 
\[
\MMM \simeq \elem(\APX(\MMM)).
\]
\end{enumerate}
\end{corollary}

\begin{proof}
Proposition \ref{targpointlem} provides the isomorphism $Emb : \MMM_\III \to \APX(\MMM)$ which proves the first part. For the second part remind that if $\MMM$ is extensional and complete as a PER-set, then $\MMM$ is maximal, extensional and complete over $\MMM_\III$, and $\MMM$ is thus a limit of $\MMM_\III$. The fact that $\elem(\MMM_\III)$ is a limit of $\MMM_\III$ and that limits are unique up to isomorphism confirms $\MMM \simeq \elem(\MMM_\III)$.
\end{proof}

Notice that each target $(\MMM,\Pmap)$ for a prefactor system $\MMM_\III$ has a unique extension $(\MMM,\Qmap)$, defined by $a \Qmap a_i \iffdef a_i \in (Ext(a))^m$, which is maximal.

\section{The Function Space of PER-Sets}
\label{funPERset}

This section is about the construction of the function space between two PER-sets $(\MMM,\per_\III, \ptp_\III)$ and $(\NNN,\per_\JJJ, \ptp_\JJJ)$. As in the function space construction for factor systems, let the index set of the function space be $\III \times \JJJ$ with product order.

\subsection{$\fil$-Continuous Functions and the Function Space of PER-Sets}
\label{funspperuniv}

We introduce now a notion of $\fil$-continuity which is similar, but not equivalent to the notion of (uniform) continuity. Given an output precision $j \in \JJJ$, then a $\fil$-continuous function $f : \MMM \to \NNN$ does not necessarily give a corresponding input precision $i \in \III$, but there is some $j' \geq j$ and some $i \in \III$ such that $a \per_i b$ implies $f(a) \per_{j'} f(b)$.

\begin{definition}[$\fil$-Continuity]
\label{sufcontdef}
Given sets with families of PERs $(\MMM,\per_\III)$ and $(\NNN,\per_\JJJ)$, then the partial equivalence relations $\per_{\III \times \JJJ}$ between two functions $f$ and $g$ from $\MMM$ to $\NNN$ are defined as a logical relation:
\begin{equation}
\label{perfuneq}
f \per_{i \to j} g \, \iffdef a \per_i b \text{ implies } f(a) \per_j g(b) \text{ for all } a,b \in \MMM.
\end{equation}

A function $f : \MMM \to \NNN$ is \emph{$i$-$j$-continuous} iff $f \in [i \to j]$ and $f$ is \emph{$\fil$-continuous}, or \emph{continuous} for short, iff there are $\fil$-many indices $i \to j$ such that $f$ is $i$-$j$-continuous.
\end{definition}

In other words, a function $f : \MMM \to \NNN$ is $i$-$j$-continuous iff $a \per_i b$ implies $f(a) \per_j f(b)$. As a consequence, $f \in [i \to j]$ and $a \in [i]$ implies $f(a) \in [j]$. From the above definition we immediately get (see Section \ref{tarsec}):

\begin{lemma}
\label{perfunoverlem}
$[\MMM \to_{\fil} \NNN] = \{f : \MMM \to \NNN \mid \III_f \in \fil(\III \times \JJJ)\}$ is the function space of all $\fil$-continuous functions. If $\MMM$ is a target for $\MMM_\III$ and $\NNN$ a target for $\NNN_\JJJ$, then $[\MMM \to_{\fil} \NNN]$ is a target for $[\MMM_\III \to \NNN_\JJJ]$.
\end{lemma}

For $\MMM$ a target for a factor system $\MMM_\III$ and $\NNN$ for $\NNN_\JJJ$, we know that $f \in [\MMM \to_\fil \NNN]$ has $\fil$-many approximations in $[\MMM_\III \to \NNN_\JJJ]$ and $a \in \MMM$ has $\fil$-many approximations in $\MMM_\III$. We also know that $f(a)$ has $\fil$-many approximations in $\NNN_\JJJ$, since $f(a) \in \NNN$ and $\NNN$ is a target for $\NNN_\JJJ$. Condition (\ref{csetcond}) makes the application of $f$ to its argument $a$ meaningful in terms of its approximations, since it guarantees that indefinitely many of these approximations of $f(a)$ stem from $\III_f$ and $\III_a$. 

\begin{example}
\label{funallex}
Consider Example \ref{natpointex} about the PER-set of the natural numbers, and let a further PER-set $\MMM$ with index set $\JJJ$ be given. Assume that $\fil(\JJJ)$ is defined in a way that it satisfies Condition (\ref{filtereq}). Define $\fil(\nat^+ \times \JJJ)$ by
\[
\HHH \in \fil(\nat^+ \times \JJJ) \iffdef \{i \in \III_\rho \mid \{j \in \JJJ \mid i \to j \in \HHH\} \in \fil(\JJJ)\} \in \fil(\nat^+),
\]
which satisfies both conditions, (\ref{csetcond}) and (\ref{filtereq}).\footnote{We will work this out in more detail in a subsequent paper, where we define these sets in general, similar as in \cite{eberl2022RML}.} 

Then every function $f: \nat \to \MMM$ is $\fil$-continuous: For all $i \in \nat^+$ we have $a \per_i b \imp a = b \imp f(a) = f(b) \imp f(a) \per_j f(b)$ for $\fil$-many $j \in \JJJ$, since $\per_\JJJ$ is dense (Definition \ref{perunivdef}). The set of all indices $i \to j$ satisfying this implication is in $\fil(\nat^+ \times \JJJ)$ by definition. In particular, each function $f : \nat \to \nat$ is $\fil$-continuous with $\{f : \nat \to \nat \mid f(n) < j \text{ for all } n < i\}$ being the $i$-$j$-continuous functions.
\end{example}

\begin{example}
\label{notcontex}
Consider $\forall \in [[\nat \to \bool] \to \bool]$ for $\nat$ as in Example \ref{natpointex} and $\bool := \{true, false\}$ with trivial index set $\{\ast\}$ and $\per_\ast$ the identity, defined as
\[
\forall : A \mapsto
\begin{cases}
true & \text{if } A(n) = true \text{ for all } n,\\
false & \text{otherwise.}
\end{cases}
\]

There is no $i \in \nat^+$ for which $\forall$ is $(i \to \ast)\,$-$\,\ast\,$-continuous: Assume that there is some $i \in \nat^+$ such that $A \per_{i \to \ast} A'$ implies $\forall(A) = \forall(A')$. But we have for the two functions $A: n \mapsto true$ and
\[
A' : n \mapsto
\begin{cases}
true & \text{if } n < i,\\
false & \text{otherwise.}
\end{cases}
\]
$A \per_{i \to \ast} A'$, as well as $\forall(A) = true \not= false = \forall(A')$. The functional $\forall$ is therefore not $\fil$-continuous. Indeed, this example shows the challenge of giving an interpretation of propositions in this model: The model uses only continuous functions, but the universal quantifier is not continuous.
\end{example}

\subsection{Points and Extensions on the Function Space}

We can easily find the points $\PPP_{i \to j}$ and the extended equivalence relation $\perp_{i \to j}$ on the function space if we consider a target for a factor system. The maps $Proj_{i \to j}$ and $Emb_{i \to j}$, defined in Section \ref{tarsec}, yield $Proj_{i \to j}(f) = Proj_j \circ f \circ Emb_i$ and $Emb_{i \to j}(f_{i \to j}) = Emb_j \circ f_{i \to j} \circ Proj_i$. So $Proj_{i \to j}(f) = Proj_{i \to j}(g)$ means that on values $a = Emb_i(a_i)$ the functions $f(a)$ and $g(a)$ are $\perp_j$-equivalent. For $f \in \PPP_{i \to j}$ this boils down to $f(a) \in \PPP_j$ and $f(a) = f(b)$ whenever $a \perp_i b$. These observations lead to the following definitions:
\begin{align*}
f \perp_{i \to j} g &\iffdef f(a) \perp_j g(a) \text{ for all } a \in \PPP_i, \\
f \in \PPP_{i \to j} &\iffdef a \perp_i b \text{ implies } f(a) = f(b) \in \PPP_j \text{ for all } a, b \in \MMM.
\end{align*}

In words, two functions $f$ and $g$ are $\perp_{i \to j}$-equivalent, if we verify on the grid points $a \in \PPP_i$ that their values $f(a)$ and $g(a)$ are $\perp_j$-equivalent. And the grid points $f \in \PPP_{i \to j}$ of the function space are the step functions with values in $\PPP_j$, constant on each of the extended equivalence classes. We immediately get the maps $\ptp_{i \to j}$:

\begin{lemma}
\label{ptlem}
The function $\ptp_{i \to j} : [\MMM \to_\fil \NNN] \to \PPP_{i \to j}$ is given by 
\begin{equation}
\label{pteq}
\ptp_{i \to j}(f) = \ptp_j \circ f \circ \ptp_i.
\end{equation}
\end{lemma}

\begin{proof}
Obviously, $\ptp_{i \to j}$ as defined in (\ref{pteq}) satisfies $\ptp_{i \to j}(f) \in \PPP_{i \to j}$ for all functions $f : \MMM \to \NNN$. Next we shall prove that $\ptp_{i \to j}(f) \perp_{i \to j} f$, i.e., it must be checked that $\ptp_j(f(\ptp_i(a))) \perp_j f(a)$ for all $a \in \PPP_i$. Clearly, $\ptp_i(a) = a$, so that $f(\ptp_i(a)) = f(a)$ and $\ptp_j(f(\ptp_i(a))) \perp_j f(a)$. 

It remains to show $\ptp_{i \to j}(f) = \ptp_{i \to j}(g)$, provided that $f \perp_{i \to j} g$. By assumption $f(\ptp_i(a)) \perp_j g(\ptp_i(a))$ holds since $\ptp_i(a) \in \PPP_i$, hence $\ptp_{i \to j}(f)(a) = \ptp_j(f(\ptp_i(a))) = \ptp_j(g(\ptp_i(a))) = \ptp_{i \to j}(g)(a)$.
\end{proof}

\begin{proposition}
$[\MMM \to_{\fil} \NNN]$ is a PER-set, provided $\MMM$ and $\NNN$ are.
\end{proposition}

\begin{proof}
We apply Proposition \ref{persetdefprop}. The properties of a PER are easily verified for $\per_{i \to j}$, and $\ptp_{i \to j}$ is surjective since $\ptp_{i \to j}(f) = f$ for all $f \in \PPP_{i \to j}$. Next observe that $\per_{\III \times \JJJ}$ is dense by the definition of $\fil$-continuous functions. The inclusions $\PPP_{i \to j} \subseteq \PPP_{i' \to j'}$ for $i \to j \leq i' \to j'$ follow from $a \perp_{i'} b \imp a \perp_{i} b \imp f(a) = f(b) \in \PPP_j \subseteq \PPP_{j'}$ for $f \in \PPP_{i \to j}$. Now we shall prove that $f \per_{i \to j} g$ implies $f \per_{i \to j} \ptp_{i \to j}(f) = \ptp_{i \to j}(g)$. 
\begin{enumerate}
\item To deduce $f \per_{i \to j} \ptp_{i \to j}(f)$, suppose $a \per_i b$. We have to show $f(a) \per_j \ptp_{j}(f(\ptp_{i}(b)))$. Clearly, $b \in [i]$, so $a \per_i b \per_i \ptp_{i}(b)$ and $f \per_{i \to j} f$ yields $f(a) \per_j f(\ptp_{i}(b)) \per_j \ptp_{j}(f(\ptp_{i}(b)))$, the latter since $f(\ptp_{i}(b)) \in [j]$.
\item We establish $\ptp_{j}(f(\ptp_{i}(a))) = \ptp_{j}(g(\ptp_{i}(a)))$ for $a \in \MMM$ by observing that $\ptp_{i}(a) \in \PPP_i \subseteq [i]$ implies $\ptp_{i}(a) \per_i \ptp_{i}(a)$. Consequently $f(\ptp_{i}(a)) \per_j g(\ptp_{i}(a))$, since $f \per_{i \to j} g$. This proves the claim.
\end{enumerate}

Let $i \to j \leq i' \to j'$. In the view of part \ref{7extperlem} of Lemma \ref{extperlem}, the forth property in Proposition \ref{persetdefprop} follows from
\[
\ptp_{i \to j}(\ptp_{i' \to j'}(f)) = \ptp_j \circ \ptp_{j'} \circ f \circ \ptp_{i'} \circ \ptp_i = \ptp_j \circ f \circ \ptp_i = \ptp_{i \to j}(f).
\]

Finally 
\begin{align*}
a \per_i b \ &\imp \ \ptp_{i'}(a) \per_i \ptp_{i'}(b) && \text{(by \ref{5extperlem}.~of Lemma \ref{extperlem})}\\
&\imp \ f(\ptp_{i'}(a)) \per_j f(\ptp_{i'}(b)) && (\text{if } f \in [i \to j]) \\
&\imp \ \ptp_{j'}(f(\ptp_{i'}(a))) \per_j \ptp_{j'}(f(\ptp_{i'}(b))) && \text{(by \ref{5extperlem}.~of Lemma \ref{extperlem})} \\
&\imp \ \ptp_{i' \to j'}(f)(a) \per_j \ptp_{i' \to j'}(f)(b) &&\text{(by Lemma \ref{ptlem}),}
\end{align*}
showing that $f \in [i \to j]$ implies $\ptp_{i' \to j'}(f) \in [i \to j]$.
\end{proof}

\subsection{Targets and Limits on the Function Space}
\label{relfunctsec}

An immediate consequence of Propositions \ref{targpointlem} and \ref{projextlem} is the following fact. 

\begin{proposition}
\label{apxisomprop}
Let $\MMM$ be a target for a factor system $\MMM_\III$ and $\NNN$ a target for a stable factor system $\NNN_\JJJ$. Then $Emb : [\MMM_\III \to \NNN_\JJJ] \to \APX([\MMM \to_\fil \NNN])$ is a strong and surjective homomorphism.
\end{proposition}

A PER-set is the target of its internal system, which automatically satisfies Property (\ref{idmapeq}) and is therefore stable. The embedding $Emb$, used in Proposition \ref{apxisomprop}, is thus an isomorphism by Proposition \ref{targpointlem}, we thus have:

\begin{corollary}
\label{apxisomlem}
For two PER-sets $\MMM$ and $\NNN$, the factor systems $\APX([\MMM \to_{\fil} \NNN])$ and $[\APX(\MMM) \to \APX(\NNN)]$ are isomorphic.
\end{corollary}

This corollary says that on the level of the factor system, the function space construction of factor systems and PER-sets are in a one-to-one correspondence. We do not need the Conditions (\ref{csetcond}) and (\ref{filtereq}) for its proof. However, both conditions are needed to prove the equivalence (\ref{elemisomeq}), i.e., to show that $[\elem(\MMM_\III) \to_{\fil} \elem(\NNN_\JJJ)]$ and $\elem([\MMM_\III \to \NNN_\JJJ])$ are also isomorphic.

Let us consider a $\fil$-continuous function $f : \MMM \to \NNN$ and assume that $\MMM$ and $\NNN$ are extensional and complete PER-sets with infinite domains $\MMM$ and $\NNN$.\footnote{This could mean that we allow actual infinite sets in our background model, so $\MMM$ and $\NNN$ are actual infinite. Otherwise, if we apply a consequent finitistic perspective, then these are potential infinite sets and we refer to $\MMM$ and $\NNN$ by referring to some indefinitely large finite sets.} So both sets are limits of factor systems $\MMM_\III$ and $\NNN_\JJJ$ resp. These could be for instance their internal factor systems (see Corollary \ref{apxlem}) or other factor systems with these limit sets. The function space $[\MMM \to_{\fil} \NNN]$ is then a limit of $[\MMM_\III \to \NNN_\JJJ]$, where a function $f \in [\MMM \to_{\fil} \NNN]$ has approximations in $[\MMM_\III \to \NNN_\JJJ]$, given by its extension $Ext(f)$, see Section \ref{tarsec}.

Suppose, for a given level of output precision $j \in \JJJ$, we want an input precision $i \in \III$. It could be that there is no $i \in \III$ such that we find an index $i \to j \in \III_f$. But we can at least find some $i \in \III$ with $i \to j' \in \III_f$ for $j' \geq j$, since $\III_f$ is in $\fil(\III \times \JJJ)$ and is thus cofinal. This provides an even better precision than required. Assume, however, that $j' = j$.

So we have an approximation $f_{i \to j} \in Ext(f)$ of $f$. The function $f$ can handle any input $a \in \MMM$, provided $i \in \III_a$. If there is no such instance $a_i$, it would be sufficient to take a more precise input $a_{i'}$ with $i' \geq i$. However, there must be an instance $f_{i' \to j}$ for it as well. Intuitively, such an instance should exist, since we allow $f$ to process a more precise input for the same level of precision as the output. For first-order functions this intuition is indeed correct and causes no problems, for higher-order functions the situation is different, see \cite{eberl2023}.

Consider again the $\fil$-continuous function $f$. Then for a given $j \in \JJJ$ we find an index $i$ such that $i \in \III_a$ and $i \to j \in \III_f$. So for a given output precision $j \in \JJJ$ and instance $a_i$ of our input $a$, we find an instance $f_{i \to j}$ of the function. The result $f(a)$ of applying $f$ to $a$ with an accuracy of $j$ can then be obtained in various ways. Since $f$ is $i$-$j$-continuous, we may also apply $f$ to the grit point next to $a$ at level $i$, that is, to $pt_i(a)$. We can also use the step function $pt_{i \to j}(f)$ instead of $f$. Definition \ref{pointperdef} and the Properties (\ref{perfuneq}) and (\ref{pteq}) show that the output is the same relative to the output precision $j \in \JJJ$:
\[
f(a) \per_j f(pt_i(a)) \per_j pt_{i \to j}(f)(a) = pt_{i \to j}(f)(pt_i(a))
\]
and $f_{i \to j}(a_i)$ is in the extension of all these values for all $a_i$ with $Emb_i(a_i) = pt_i(a)$, and for all $f_{i \to j}$ such that $Emb_{i \to j}(f_{i \to j}) = pt_{i \to j}(f)$.

\section{Conclusion and Future Work}

In this paper, we introduce PER-sets as spaces consisting of regions of arbitrarily small size, as opposed to points, as their fundamental components. The size of these regions depends on the meta-level context of the mathematical investigation. As mentioned in Section \ref{observesec}, proving continuity only requires this indefinite partition of space, rather than an infinite one. A future objective is to demonstrate that all properties are satisfied by this indefinite partition. To achieve this, type theory will serve as a framework. Therein the type $bool$ represents propositions, which are interpreted as the two truth values. Since higher-order logic can be embedded into simple type theory, this guarantees that all expressible propositions have a potential infinite reading, not only first-order propositions, but also higher-order propositions. This integration of logic will be presented in a subsequent paper, where a general reflection principle for a fragment of (simple and later also dependent) type theory will be shown. A first step into this direction has been done in \cite{eberl2024}.

\appendix
\section{Postponed Proofs}

Here we present proofs that we postponed for the sake of readability.

\subsection{Proofs of Section \ref{richuni} - \nameref{richuni}}

\begin{proofof}{Proposition \ref{targpointlem}}

$\MMM$ is a $\fil$-set by Lemma \ref {approxlem}. Before showing that it is pointed, we shall prove the additional claims: The $\imp$ part of the first claim holds since the $\comp$-embeddings $Emb_i$ preserve $\comp$. Conversely, $a_i \comp b_i$ is a consequence of (\ref{ppmapeq}) applied to $Emb_i(a_i) \Pmap a_i$ and $Emb_i(a_i) = Emb_i(b_i) \Pmap b_i$. 

For the second claim let $a_{i'} \pmap a_i$, then $Emb_{i'}(a_{i'}) \Pmap a_i$ follows by (\ref{embtareq}) and hence $pt_i(Emb_{i'}(a_{i'})) = Emb_i(a_i)$ by (\ref{perpmap2}). Conversely, let $\MMM_\III$ be stable and assume $pt_i(Emb_{i'}(a_{i'})) = Emb_i(a_i)$. It must be checked that $a_{i'} \pmap a_i$. Because of (\ref{perpmap2}) we can find some $b_i \in \MMM_i$ such that $Emb_{i'}(a_{i'}) \Pmap b_i$ and $Emb_i(a_i) = Emb_i(b_i)$. Now observe that $Emb_{i'}(a_{i'}) \Pmap a_{i'}$ by (\ref{embcond}) so that $a_{i'} \pmap b_i$ follows from (\ref{ppmapeq}). From (\ref{perpmapcons}) and $Emb_i(a_i) = Emb_i(b_i)$ we conclude $a_i \comp b_i$ and because $\MMM_\III$ is stable we have $a_{i'} \pmap a_i$. 

With these additional claims we are able to show that the map $pt_i : [i] \to \PPP_i$ is well defined and surjective. If $a \Pmap a_i$ and $a \Pmap b_i$, then $a_i \comp b_i$ follows from (\ref{perpmap}). Now $Emb_i(a_i) = Emb_i(b_i)$ is a consequence of the first claim of this proposition, and $pt_i$ is surjective by the second claim, applied to $a_i \pmap a_i$. Next we have to show the properties stated in Definition \ref{pointperdef}. 

\begin{enumerate}
\item If $a \in \PPP_i$ then $a = Emb_i(a_i)$ for some $a_i \in \MMM_i$ and consequently $a = Emb_i(a_i) = Emb_{i'}(\emb{i}{i'}(a_i)) \in \PPP_{i'}$.

\item If $a \per_i b$, then $a \Pmap a_i$ and $b \Pmap b_i$ for some $a_i \comp b_i$ in $\MMM_i$ by (\ref{perpmap}). So $pt_i(a) = Emb_i(a_i) \Pmap a_i$, which shows $a \per_i pt_i(a)$. Moreover, $a_i \comp b_i$ implies $Emb_i(a_i) = Emb_i(b_i)$, i.e., $pt_i(a) = pt_i(b)$.

\item Assume $a \in [i'] \cap [i]$, then there are $a_{i'} \in \MMM_{i'}$ and $a_i \in \MMM_i$ such that $a \Pmap a_{i'}$ and $a \Pmap a_i$, hence $a_{i'} \pmap a_i$ by Property (\ref{pmapeq}). From $a_{i'} \pmap a_i$ and (\ref{embtareq}) we conclude $Emb_{i'}(a_{i'}) \Pmap a_i$. Furthermore, we have $pt_{i'}(a) = Emb_{i'}(a_{i'})$, because $a \Pmap a_{i'}$. Consequently $pt_{i'}(a) \Pmap a_i \in \MMM_i$, showing that $pt_{i'}(a) \in [i]$.
\end{enumerate}

For the final statement assume that $\MMM_\III$ is stable and remind that $\APX(\MMM)$ consists of the points $\PPP_\III$. Each $Emb_i : \MMM_i \to \PPP_i$ is obviously surjective. Because $\MMM_\III$ is stable, $Emb$ preserves $\pmap$ in both directions by the second claim, and $Emb$ preserves embeddings since $Emb_{i'}(\emb{i}{i'}(a_i)) = Emb_i(a_i)$. In case (\ref{idmapeq}) holds, $Emb$ is injective by the first claim, and hence an isomorphism.

\end{proofof}

\subsection{Proofs of Section \ref{eqrelpersec} - \nameref{eqrelpersec}}

\begin{proofof}{Proposition \ref{projextlem}}

In addition to the proof of Proposition \ref{targpointlem}, the properties from Definition \ref{extperdef} must be checked.

\begin{enumerate}
\item First we have to show that $a \Pmap a_i \in \MMM_i$ implies $Emb_i(Proj_i(a)) = Emb_i(a_i)$. From Section \ref{tarsec} we use the fact that if $a \Pmap a_i$, then $Proj_i(a) \pmap a_i$, which is the same as $Proj_i(a) \comp a_i$, and hence we conclude that $Emb_i(Proj_i(a)) = Emb_i(a_i)$. The relation $Proj_i(a) \comp a_i$ immediately shows the additional claim (\ref{projembeq}) as well, if we apply it to $a := Emb_{i}(a_{i}) \Pmap a_{i}$.

\item To prove $\ptp_i = \ptp_i \circ \ptp_{i'}$ for $i' \geq i$ means to show 
\[
Emb_i(Proj_i(a)) = Emb_i(Proj_i(Emb_{i'}(Proj_{i'}(a)))).
\]

Let $a_{i'} := Proj_{i'}(a)$, so that $Proj_{i'}(a) = a_{i'} \comp Proj_{i'}(Emb_{i'}(a_{i'}))$. By the definition of relation $\perp_{i'}$ this is  $a \perp_{i'} Emb_{i'}(a_{i'})$, which implies $a \perp_i Emb_{i'}(a_{i'})$ by Lemma \ref{extperlem}. Again, by the definition of $\perp_i$, we infer $Emb_i(Proj_i(a)) = Emb_i(Proj_i(Emb_{i'}(a_{i'})))$.

\item To confirm that $a \in [i]$ implies $\ptp_{i'}(a) \in [i]$, suppose $a \Pmap a_i \in \MMM_i$. Then $Proj_{i'}(a) \pmap a_i$ and $Emb_{i'}(Proj_{i'}(a)) \Pmap a_i$ by (\ref{embtareq}). This establishes $\ptp_{i'}(a) \in [i]$.
\end{enumerate}

It remains to show that $Emb_\III$ is a factor system homomorphism, which means that it additionally commutes with $proj_\III$. Indeed, we have
\[
Emb_i(\proj{i'}{i}(a_{i'})) = Emb_i(Proj_i(Emb_{i'}(a_{i'}))) = \ptp_i(Emb_{i'}(a_{i'})).
\]

The first equation follows from (\ref{Projprop}) and the fact that $Emb$ preserves $\comp$, the second equation is the definition of $\ptp_i$.
\end{proofof}

\bibliographystyle{elsarticle-num}
\bibliography{ContPotInfMod}

\end{document}
\endinput